\RequirePackage{fix-cm}

\documentclass[smallextended]{svjour3opt} 
\usepackage{etex}
\usepackage{todonotes}
\usepackage{lscape} 
\usepackage{palatino,url}
\usepackage{verbatim, amssymb,amsmath}
\usepackage{varwidth}
\usepackage[export]{adjustbox}
\usepackage{multirow}

\usepackage{mathrsfs}
\usepackage{graphics,amsopn,amstext,amsfonts,color}
\usepackage{bm}
\usepackage{setspace}
\usepackage{subfig}
\usepackage{graphicx}
\usepackage{epstopdf}
\usepackage{caption}
\usepackage{multirow}
\usepackage{booktabs}
\usepackage{paralist}
\usepackage{colortbl}
\usepackage{appendix}
\usepackage{soul,caption,subfig}
\usepackage[bookmarksnumbered=true]{hyperref}
\usepackage{graphics,latexsym,amsfonts,verbatim,amsmath}
\usepackage{algorithm}
\usepackage{algorithmicx}
\usepackage{algpseudocode}
\usepackage{tikz}
\usetikzlibrary{calc}
\usetikzlibrary{plotmarks}

\usepackage{pgfplots}

\usepackage{pgf}
\usepackage{pst-func}

\definecolor{pinegreen}{rgb}{0.0, 0.47, 0.44}
\usepackage{pgfpages}
\usepackage{todonotes}
\usepackage{hhline} 
\usepackage{empheq}
\usepackage{cases}

\def\L{{\mathcal L}}

\def\E{{\mathbb E}}
\def\Var{{\mathbb V}}

\def\Pr{{\mathbb{P}}}

\def\Re{\mathbb{R}}
\def\Qe{\mathbb{Q}}

\def\hat{\widehat}

\def \bxi{\boldsymbol{\xi}}
\def \bzeta{\boldsymbol{\zeta}}

\def \P{\mathcal{P}}

\def \Z{{\mathcal{Z}}}

\def\L{{\mathcal L}}

\def\N{{\mathcal N}}

\def\P{{\mathcal P}}

\def\Re{{\mathbb R}}

\newcommand{\exclude}[1]{}

 \newcommand{\y}{\mathbf{y}}

\newcommand{\e}{\mathbf{e}}

\renewcommand{\proof}{\noindent{\em Proof: }}

\newcommand*{\QEDA}{\hfill\ensuremath{\square}}
\newcommand*{\QEDB}{\hfill\ensuremath{\diamond}}

\def\delequal{\stackrel{\triangle}{=}} 

\algdef{SE}[DOWHILE]{Do}{doWhile}{\algorithmicdo}[1]{\algorithmicwhile\ #1}%

\def\L{{\mathcal L}}

\def\E{{\mathbb E}}
\def\Var{{\mathbb V}}

\def\Pr{{\mathbb{P}}}

\def\Re{\mathbb{R}}
\def\Qe{\mathbb{Q}}

\def\hat{\widehat}

\def \P{\mathcal{P}}

\def \Z{{\mathcal{Z}}}

\def\L{{\mathcal L}}

\def\N{{\mathcal N}}

\def\P{{\mathcal P}}

\def\Re{{\mathbb R}}

\def\bxi{\tilde{\boldsymbol \xi}}

\def\bfxi{\boldsymbol \xi}
\def\bzeta{\boldsymbol \zeta}
\def\x{\vect{x}}

\newcommand{\vect}[1]{\boldsymbol{\bm{#1}}}

\newcommand{\subparagraph}{}
\usepackage{titlesec}

\newtheorem{assumption}{Assumption}

\usepackage[numbers,sort&compress,comma]{natbib}
\allowdisplaybreaks

\title{Distributionally Robust Bottleneck Combinatorial Problems: Uncertainty Quantification and Robust Decision Making}
\date{\today}
\titlerunning{Distributionally Robust Bottleneck Combinatorial Problems}

\author{Weijun Xie \and Jie Zhang \and Shabbir Ahmed}

\institute{Weijun Xie \at
	Virginia Tech, Blacksburg, VA\\
	\email{wxie@vt.edu}
	\and Jie Zhang \at
	Virginia Tech, Blacksburg, VA\\
	\email{jiezhang@vt.edu}
	\and 
	Shabbir Ahmed \at
	Georgia Institute of Technology, Atlanta, GA\\
	\email{sahmed@isye.gatech.edu}
}

\allowdisplaybreaks

\makeatletter
\pgfmathdeclarefunction{erf}{1}{%
	\begingroup
	\pgfmathparse{#1 > 0 ? 1 : -1}%
	\edef\sign{\pgfmathresult}%
	\pgfmathparse{abs(#1)}%
	\edef\x{\pgfmathresult}%
	\pgfmathparse{1/(1+0.3275911*\x)}%
	\edef\t{\pgfmathresult}%
	\pgfmathparse{%
		1 - (((((1.061405429*\t -1.453152027)*\t) + 1.421413741)*\t
		-0.284496736)*\t + 0.254829592)*\t*exp(-(\x*\x))}%
	\edef\y{\pgfmathresult}%
	\pgfmathparse{(\sign)*\y}%
	\pgfmath@smuggleone\pgfmathresult%
	\endgroup
}
\makeatother

\begin{document}
	\maketitle

	\begin{abstract} 
		In a bottleneck combinatorial problem, the objective is to minimize the highest cost of elements of a subset selected from the combinatorial solution space. This paper studies data-driven distributionally robust bottleneck combinatorial problems (DRBCP) with stochastic costs, where the probability distribution of the cost vector is contained in a ball of distributions centered at the empirical distribution specified by the Wasserstein distance. We study two distinct versions of DRBCP from different applications: (i) Motivated by the multi-hop wireless network application, we first study the uncertainty quantification of DRBCP (denoted by DRBCP-U), where decision-makers would like to have an accurate estimation of the worst-case value of DRBCP. The difficulty of DRBCP-U is to handle its max-min-max form. Fortunately, similar to the strong duality of linear programming, the alternative forms of the bottleneck combinatorial problems using clutters and blocking systems allow us to derive equivalent deterministic reformulations, which can be computed via mixed-integer programs. In addition, by drawing the connection between DRBCP-U and its sampling average approximation counterpart under empirical distribution, we show that the Wasserstein radius can be chosen in the order of negative square root of sample size, improving the existing known results; and (ii) Next, motivated by the ride-sharing application, decision-makers choose the best service-and-passenger matching that minimizes the unfairness. That is, we study the decision-making DRBCP, denoted by DRBCP-D. For DRBCP-D, we show that its optimal solution is also optimal to its sampling average approximation counterpart, and the Wasserstein radius can be chosen in a similar order as DRBCP-U. When the sample size is small, we propose to use the optimal value of DRBCP-D to construct an indifferent solution space and propose an alternative decision-robust model, which finds the best indifferent solution to minimize the empirical variance. We further show that the decision robust model can be recast as a mixed-integer conic program. Finally, we extend the proposed models and solution approaches to the distributionally robust $\Gamma-$sum bottleneck combinatorial problem (DR$\Gamma$BCP), where decision-makers are interested in minimizing the worst-case sum of $\Gamma$ highest costs of elements.
	\end{abstract}

	\newpage
	\section{Introduction}
	\subsection{Setting}
	
	In a combinatorial problem, there are $n$ elements, indexed by $[n]:=\{1,\ldots, n\}$, and a collection of subsets $X\subseteq 2^{[n]}$, where each element $j\in [n]$ bears a cost $c_j$. The objective of a bottleneck combinatorial problem (BCP) \cite{kasperski2013bottleneck} is to seek the best subset $x\in X$ that minimizes its highest cost of elements (i.e., bottleneck cost), defined as $\min_{x\in X}\max_{j\in x}c_j$, i.e., a BCP can be formulated as
	\begin{align}
	\textrm{(BCP)}\quad Z(\bm{c})=\min_{x\in X}\max_{j\in x}c_j. \label{eq_cbp}
	\end{align}
	The BCP \eqref{eq_cbp} covers many interesting bottleneck problems, e.g., bottleneck shortest path \cite{kaibel2006bottleneck}, bottleneck spanning tree \cite{camerini1978min}, bottleneck assignment \cite{spivey2011asymptotic}, etc.

	Oftentimes, the cost vector is random, denoted by $\tilde{\bm{c}}$, and its probability distribution $\Pr$ is also difficult to characterize. Therefore, this paper considers the Distributionally Robust Combinatorial Bottleneck Problem (DRBCP). Let ambiguity set $\P$ represent a distributional family. In this paper, we focus on two interesting variants of DRBCP: uncertainty quantification (DRBCP-U) and decision-making (DRBCP-D), i.e.,
	\begin{subequations}
		\begin{align}
		\textrm{(DRBCP-U)}\quad&v_{U}=\sup_{\Pr\in \P}\E_{\Pr}\left[Z(\tilde{\bm{c}})\right],\quad \textrm{and}\label{general_eq_dis_robust}\\
		\textrm{(DRBCP-D)}\quad &v_{D}=\min_{x\in X}\sup_{\Pr\in \P}\E_{\Pr}\left[\max_{j\in x}\tilde c_j\right].\label{general_eq_dis_robust_d}
		\end{align}
	\end{subequations}
	For DRBCP-U \eqref{general_eq_dis_robust}, the objective is to determine the value of the worst-case bottleneck cost, while in DRBCP-D \eqref{general_eq_dis_robust_d}, it aims to find the best subset to minimize the worst-case bottleneck cost. Both DRBCP-U \eqref{general_eq_dis_robust} and DRBCP-D \eqref{general_eq_dis_robust_d} are motivated from real-world applications. For example, in a wireless sensor network (e.g., multi-hop network) \cite{liu2005multi, gao2006determining, yoo2007maximum, zhang2008achieving}, there is usually a pair of transmitter (Tx) and receiver (Rx) and throughout the communication network, there is a minimum data rate requirement to be strictly enforced. Thus, it is important for the wireless network coordinator to determine the data rate between Tx and Rx, i.e., solve DRBCP-U \eqref{general_eq_dis_robust}. Another example is that in the bottleneck matching problem \cite{spivey2011asymptotic} (e.g., online matching between ride-sharing service providers and passengers), since the travel time is often stochastic, it is important to find a fair matching between ride-sharing service providers and passengers which minimizes the longest waiting time of passengers, i.e., solve DRBCP-D \eqref{general_eq_dis_robust_d}. Besides, DRBCP-D \eqref{general_eq_dis_robust_d} can be also applied to organ transplanting \cite{parent2017fair} to fairly distribute available organs to the patients, or to nurses' scheduling problems \cite{martin2013cooperative} to minimize the longest overtime workloads.

	\subsection{Wasserstein Ambiguity Set}
	Given an empirical distribution $\Pr_{\tilde{\bzeta}}$ constructed using i.i.d. historical data $\Z=\{\bar{\bm{c}}^k\}_{k\in [N]}$ such that $\Pr_{\tilde{\bzeta}}\{\tilde{\bzeta}=\bar{\bm{c}}^k\}=1/N$, this paper considers data-driven Wasserstein ambiguity set (see, e.g., \cite{gao2016distributionally,blanchet2019confidence,esfahani2018data,blanchet2019quantifying,hanasusanto2018conic,chen2018data,xie2018distributionally,abadeh2018wasserstein,kuhn2019wasserstein,blanchet2019robust,chen2019sharing} as below:
	\begin{align}\label{eq_general_das}
	\P_q=\left\{\Pr: W_q\left(\Pr,\Pr_{\tilde{\bzeta}}\right)\leq \theta\right\},
	\end{align}
	where $\theta\geq 0$ denotes the Wasserstein radius and for any $q\in [1,\infty]$, the Wasserstein distance $W_q(\cdot,\cdot)$ is defined as
	\[W_q\left(\Pr_1,\Pr_{2}\right)=\inf\left\{\sqrt[q]{\int_{\Xi\times\Xi}\|{\bm{\xi}}_1-{\bm{\xi}}_2\|^q \Qe(d\bm{\xi}_1,d\bm{\xi}_2)}:\begin{array}{l}\text{$\Qe$ is a joint distribution of $\tilde{\bm{\xi}}_1$ and $\tilde{\bm{\xi}}_2$}\\
	\text{with marginals $\Pr_1$ and $\Pr_2$, respectively}\end{array}\right\}.\]
	In particular, when $q=\infty$, it reduces to the $\infty-$Wasserstein distance
	\[W_\infty\left(\Pr_1,\Pr_{2}\right)=\inf\left\{Q-\text{ess} \sup\|\tilde{\bm\xi}_1-\tilde{\bm\xi}_2\|:\begin{array}{l}\text{$\Qe$ is a joint distribution of $\tilde{\bm{\xi}}_1$ and $\tilde{\bm{\xi}}_2$}\\
	\text{with marginals $\Pr_1$ and $\Pr_2$, respectively}\end{array}\right\}.\]
	Above, $Q$-ess sup $ \left \| \cdot  \right \|$ is the essential supreme of $ \left \| \cdot  \right \|$ with respect to the joint distribution $Q$, which is formally defined as 
	$$Q-\text{ess} \sup\|\tilde{\xi}_1-\tilde{\xi}_2\|\delequal \inf \left\{\Delta: Q[\|\tilde{\xi}_1-\tilde{\xi}_2\|>\Delta]=0\right\}.$$
	It has been shown in \cite{bertsimas2018data} that its corresponding $\infty-$ Wasserstein ambiguity set $\P_\infty$ has the following equivalent form
	\begin{align}
	\P_\infty=\left\{\frac{1}{N}\sum_{k\in [N]}\delta(\tilde{\bm c}-\bm{c}^k):\exists \bm{c}^k, \|\bm{c}^k-\bar{\bm{c}}^k\|\leq \theta\right\},\label{eq_dual_wass}
	\end{align}
	where $\delta(\cdot)$ is the Dirac delta function. This neat representation has a straightforward interpretation, i.e., the worst-case distribution is also supported by $N$ points and each support point can only deviate at most $\theta$ amount from one of the empirical data $\Z=\{\bar{\bm{c}}^k\}_{k\in [N]}$. 
	
	Since $W_\infty(\cdot,\cdot)$ has the similar statistical performance as $W_q(\cdot,\cdot)$ with $q\in [1,\infty)$ \cite{fournier2015rate,trillos2015rate}, we will focus on $\infty-$Wasserstein ambiguity set, i.e., throughout this paper, we make the following assumption:
	\begin{assumption}Suppose that $\P=\P_{\infty}$.
	\end{assumption}
	The main results derived in this paper can be also extended to the other $q-$Wasserstein ambiguity set with $q\in [1,\infty)$, which is shown in the appendix. The advantages of using $\infty-$Wasserstein ambiguity set can be found in recent works (see, e.g., \cite{gao2016distributionally,bertsimas2018data,xie2019tractable}). Particularly, this paper will show that for DRBCP-U \eqref{general_eq_dis_robust}, the resulting formulation using $\infty-$Wasserstein ambiguity set can be decomposed into $N$ separate mixed-integer programs, and thus can be computed very effectively, while for the other $q-$Wasserstein ambiguity set with $q\in [1,\infty)$, the decomposition is not straightforward and requires additional computational efforts.

	\subsection{Relevant Literature}
	
	Bottleneck combinatorial problems (BCPs) have been widely studied in literature \cite{kasperski2013bottleneck,kaibel2006bottleneck,camerini1978min,spivey2011asymptotic,shinn2013efficient}. Existing approaches mainly focused on BCPs with deterministic costs and developing polynomial-time algorithms. Recently, some development has been made to study BCPs with random costs, which are mainly focused on algorithms and approximations of some special variants of stochastic BCPs. 
	We separate these works according to DRBCP-U \eqref{general_eq_dis_robust} and DRCBP-D \eqref{general_eq_dis_robust_d}, respectively.
	\begin{itemize}[(i)]
		\item[Uncertainty Quantification, i.e., DRBCP-U \eqref{general_eq_dis_robust}.] 
		Existing works attempted to derive the asymptotic behavior of stochastic BCP (see, e.g., \cite{pferschy1996random}), where the probability distribution of the random costs is given, and the number of elements $n$ tends to be infinite. In \cite{spivey2011asymptotic}, the author derived all of the asymptotic moments of the stochastic bottleneck assignment problem in which edge costs are i.i.d. from a continuous distribution. Similar results can be found in \cite{pferschy1996random}, where the author assumed the distribution is uniform. The work \cite{albrecher2005note} investigated the asymptotic gap between the minimum pessimistic bottleneck cost (i.e, $\min_{x\in X}\max_{j\in x}\tilde{c}_j$) and the minimum optimistic bottleneck cost (i.e, $\min_{x\in X}\min_{j\in x}\tilde{c}_j$) under the probabilistic setting. All of these works assumed that the probability distribution is known, and the random cost of each element is i.i.d., which might not be realistic. In this paper, we do not assume known probability distribution, and instead, decision-makers can obtain empirical data to construct an ambiguity set. 
			\item[Decision Making, i.e., DRBCP-D \eqref{general_eq_dis_robust_d}.] Many existing works focused on stochastic counterpart of DRBCP-D \eqref{general_eq_dis_robust_d}, where a probability distribution is assumed to be known. For example, in \cite{kasperski2013bottleneck}, the authors proposed new criteria to generalize stochastic BCPs with discrete random costs. The complexity analysis and approximation results for several special cases were also provided. A similar result can be found in \cite{kasperski2011possibilistic}. In \cite{perez2018algorithms}, the authors developed algorithms to solve minimax regret bottleneck path problems with interval uncertainty. The work \cite{stancu1999stochastic} studied risk-averse BCPs and reduced them to deterministic BCPs. In \cite{shen2015chance}, the authors formulated the stochastic bottleneck spanning trees problem with independently distributed edge weights as a chance constrained program and provided efficient solution algorithms. All of the works mentioned above assumed the known probability distribution, while this paper studies DRBCP-D \eqref{general_eq_dis_robust_d} with unknown probability distributions. As far as we are concerned, only in \cite{natarajan2013probabilistic}, the authors studied a generalization of DRBCP-D to minimize the worst-case sum of $\Gamma$ highest costs of elements with known marginal distributions, while different from \cite{natarajan2013probabilistic}, this paper studies both DRBCP-U and DRBCP-D under the less conservative Wasserstein ambiguity set. 
		
	\end{itemize}

	As an effective way for the decision-making under uncertainty without fully knowing the probability distribution, distributionally robust optimization has attracted much attention \cite{delage2010distributionally,Xie2016drccp,zymler2013distributionally,gao2016distributionally,hanasusanto2015distributionally,hanasusanto2015Ambiguous,blanchet2019confidence,chen2018data,xie2018distributionally}. Interested readers are referred to \cite{rahimian2019DROreview} for a comprehensive review. However, existing works mainly studied deterministic reformulations of distributionally robust optimization, where the objective function is in the form of either piece-wise convex, piece-wise concave or two-stage recourse function. This paper presents two unique aspects for distributionally robust optimization: (i) for DRBCP-U \eqref{general_eq_dis_robust}, existing deterministic reformulation techniques are not directly applicable due to the complicated minimax operator and thus we develop new techniques, allowing us to represent DRBCP-U \eqref{general_eq_dis_robust} as mixed-integer program (MIP); and (ii) for DRBCP-D \eqref{general_eq_dis_robust_d}, we observe that the sampling average approximation of a stochastic BCP is optimal, which, however, might not generate robust decisions. Thus, we present a better alternative for DRBCP-D \eqref{general_eq_dis_robust_d}, which simultaneously guarantees the robustness of decision as well as the conservatism of objective value.

	\subsection{Summary of Contributions} 
	Since in general solving a BCP is NP-hard \cite{hsu1979easy}, this paper will focus on the deterministic reformulations of DRBCP-U \eqref{general_eq_dis_robust} and DRBCP-D \eqref{general_eq_dis_robust_d} under $\infty-$Wasserstein ambiguity set $\P_{\infty}$ and derive their mixed integer-programming representations. Our main contributions are summarized as below:
	\begin{enumerate}[(i)]
		\item By refining set $X$ as an equivalent collection of mutually noncomparable subsets, we are able to derive an equivalent form of BCP \eqref{eq_cbp}. Using the equivalent representations of BCP \eqref{eq_cbp} and $\P_\infty$, we show that DRBCP-U \eqref{general_eq_dis_robust} can be formulated as $N$ separate MIPs. 
		\item By drawing the connection between DRBCP-U \eqref{general_eq_dis_robust} and its sampling average approximation counterpart, we show that the Wasserstein radius can be chosen in the order of $O(N^{-1/2})$, i.e., $\theta=O(N^{-1/2})$. This observation demonstrates that one does not need too many samples to ensure that with high probability, the optimal value of DRBCP-U \eqref{general_eq_dis_robust} provides an upper bound to the true stochastic BCP.
		\item We observe that the optimal decision of DRBCP-D \eqref{general_eq_dis_robust_d} is always equal to its sampling average approximation counterpart. Thus, we propose a new decision-robust model, which simultaneously guarantees the robustness of the objective value and the robustness of the decision.
		\item We show that the proposed results can be extended to the distributionally robust $\Gamma-$sum bottleneck combinatorial problem (DR$\Gamma$BCP), whose objective is to minimize the worst-case sum of $\Gamma-$ highest costs of elements.
	\end{enumerate}
	
	\noindent\textbf{Notation.}  The following notation is used throughout the paper. We use bold-letters (e.g., $\vect{x},\vect{A}$) to denote vectors and matrices and use corresponding non-bold letters to denote their components. Given a vector or matrix $\bm{x}$, its zero norm $\left \| \bm{x}\right \|_0$ denotes the number of its nonzero elements. We let $\e$ be the vector or matrix of all ones, and let $\e_i$ be the $i$th standard basis vector.
	Given an integer $n$, we let $[n]:=\{1,2,\ldots,n\}$, and use $\Re_+^n:=\{\vect{x}\in \Re^n:x_i\geq0, \forall i\in [n]\}$. Given a real number $t$, we let $(t)_+:=\max\{t,0\}$. Given a finite set $I$, we let $|I|$ denote its cardinality. We let $\bxi$ denote a random vector and denote its realizations by $\bfxi$.  Given a probability distribution $\Pr$ on $\Xi$, we use $\Pr\{A\}$ to denote $\Pr\{\bxi: \text{condition} \ A(\bxi) \ \text{holds}\}$ when $A(\bxi)$ is a condition on $\bxi$, and to denote $\Pr\{\bxi: \bxi \in A\}$ when $A \subseteq \Xi$. Given a set $S$ and a positive integer $\Gamma$, we let ${S\choose\Gamma}$ denote the collection of all the size-$\Gamma$ subsets of $S$. We use superscript $k\in [N]$ denote the index of scenario $k$. We let $\delta(\cdot)$ denote a Dirac delta function, roughly speaking, $\delta(t)=\begin{cases}
	\infty,&\text{if }t=0,\\
	0,&\text{otherwise}.
	\end{cases}$ 
	Additional notation will be introduced as needed.\\

	\noindent\textbf{Organization.}
	The remainder of the paper is organized as follows. In Section \ref{example}, we introduce the definition of clutters and blocking systems and an equivalent formulation of BCP. Sections \ref{sec_DRBCP-U} and \ref{sec_DRBCP-D} show the equivalent reformulations and properties of RBCP-U \eqref{general_eq_dis_robust} and DRBCP-D \eqref{general_eq_dis_robust_d} under $\infty-$Wasserstein ambiguity set $\P_{\infty}$. Section \ref{sec_numerical} shows the numerical illustration and Section \ref{sec_T_sum_dro} extends the results to distributionally robust $\Gamma$-sum bottleneck problems. Finally, Section~\ref{sec_conclusion} concludes this paper.

	\section{Clutter and Equivalent Dual Representation of BCP \eqref{eq_cbp}}\label{example}
	In this section, we will introduce important concepts and preliminary results useful for the main derivations in the following sections.
	
	\subsection{Clutter and Blocking System}
	Following from \cite{edmonds1970bottleneck,hsu1979easy,schrijver2003combinatorial}, we first define clutters and blocking systems as follows.
	\begin{definition}
		(Clutter)  A clutter $\overline{X}\subseteq 2^{[n]}$ on the set $[n]$ is a family of mutually noncomparable subsets of $[n]$, i.e., 
		$$\overline{X}=\{h \subseteq [n]: \ \text{there is no subset} \ q\in \overline{X} \ \text{such that}\  q \subsetneq h\}.$$
	\end{definition}
	
	
	\begin{definition}
		(Blocking System) Given a clutter $\overline{X}$ on $[n]$, suppose that the unique set $F$ consists of the minimal subsets of $[n]$ that have non-empty intersection with every member of $\overline{X}$, i.e.,
		$$F=\{y\subseteq [n]: y\cap h\neq \emptyset, \forall h\in\overline{X} \ \text{and}\ y \ \text{is minimal with this property}\}.$$
		Then the pair $(\overline{X}, F)$ is called a \textit{blocking system} on $[n]$ and $F=b(\overline{X})$ is called the \text{blocker} of $\overline{X}$.
	\end{definition}
	According to \cite{edmonds1970bottleneck}, for any clutter $\overline{X}$ on $[n]$, its blocker always exists and is unique.
	
	
	The following examples from \cite{edmonds1970bottleneck,hsu1979easy} list some commonly-used blocking systems.
	\begin{example}\label{example1}
		(Blocking System of the Cardinality Set) Consider a collection of size-$m$ subsets $\overline{X}=\{x\subseteq [n]: |x|=m\}$, then its blocker $F=b(\overline{X})=\{y\subseteq [n]: |y|=n-m+1\}$. 
	\end{example}
	\begin{example}\label{example2}
		(Blocking System of $s-t$ Paths) Suppose the clutter $\overline{X}$ is a collection of all the simple paths from vertex $s$ to vertex $t$ in a graph $G(V,E)$, then its blocker $F=b(\overline{X})$ is the collection of all the $s-t$ cuts separating vertices $s$ and $t$ in the graph. Here, a $s-t$ cut is defined as a collection of edges between a partition of nodes $V$, i.e., two subsets of nodes $V_1,V_2$ such that $V_1\cap V_2=\emptyset, V_1\cup V_2=V$, and $s\in V_1, t\in V_2$.
	\end{example}
	
	\begin{example}\label{example3}
		(Blocking System of Spanning Trees) Suppose the clutter $\overline{X}$ is a collection of all the spanning trees in an undirected graph $G(V,E)$, then its blocker $F=b(\overline{X})$ is the collection of all the cuts in the graph. Here, a cut is defined as a collection of edges between a partition of nodes $V$, i.e., two subsets of nodes $V_1,V_2$ such that $V_1\cap V_2=\emptyset, V_1\cup V_2=V$.
	\end{example}
	
	\begin{example}\label{example4}
		(Blocking System of the Assignments) Suppose the clutter $\overline{X}$ is a collection of all the possible $[m]\times [m]$ assignments, then its blocker $F=b(\overline{X})=\{\text{All the $h\times k$ submatrices}: h,k\subseteq [m],|h|+|k|=m+1\}$. 
	\end{example}

	\subsection{Equivalent Representation of BCP \eqref{eq_cbp}}\label{sec_dual_bcp}
	To derive an equivalent representation of BCP \eqref{eq_cbp}, let us first define a set $\overline{X}=X\setminus \{x\in X: \exists y\in X, y\subsetneq x\}$, i.e., $\overline{X}$ is the minimal subset of $X$ that each member of $X$ is either a subset or equal to a member of $\overline{X}$ but no member of $\overline{X}$ is strictly contained in another member of $\overline{X}$. Clearly, we have
	\begin{lemma}\label{cluter_blocker_-0}
		Set $\overline{X}$ is a clutter. 
	\end{lemma}
	\proof For any $x\in \overline{X}$, there is no $y\in \overline{X}$ such that $y \subsetneq x$. Thus, $\bar{X}$ is a clutter. \QEDA

	Further, in BCP \eqref{eq_cbp}, we can achieve the same optimal value if we replace $X$ by $\bar{X}$.
	\begin{lemma}\label{cluter_blocker_0}BCP \eqref{eq_cbp} is equivalent to
		\begin{equation} \label{equa_clust}
		Z(\bm{c})=\min_{h\in \overline{X}}\max_{i\in h}c_i.
		\end{equation}
	\end{lemma}
	\proof Let us denote $\hat{Z}(\bm{c})=\min_{h\in \overline{X}}\max_{i\in h}c_i$. Since $ \overline{X}\subseteq X$, we must have $\hat{Z}(\bm{c})\geq Z(\bm{c})$.
	
	On the other hand, according to the definition of $\overline{X}$, for each $x\in X $, there must exist an $h\in \overline{X}$ such that $h\subseteq x$, which implies that $\max_{i\in h}c_i \leq \max_{i\in x}c_i$. Thus, we must have $\min_{h\in \overline{X}}\max_{i\in h}c_i \leq \min_{x\in X}\max_{i\in x}c_i$, i.e., $\hat{Z}(\bm{c})\leq Z(\bm{c})$.
	\QEDA
	
	Finally, we present the equivalent formulation of \eqref{equa_clust}.
	\begin{theorem}\label{cluter_blocker}Suppose $F$ is the blocker of $\overline{X}$. Then BCP \eqref{eq_cbp} is equivalent to
		\begin{equation}\label{block}
		Z(\bm{c})=\max_{y\in F}\min_{j\in y}c_j.
		\end{equation}
	\end{theorem}
	\proof According to Lemma~\ref{cluter_blocker_-0}, $\overline{X}$ is a clutter and its blocker $F$ exists and is unique. Then the conclusion follows from Lemma~\ref{cluter_blocker_0} and the theorem on page 300 in \cite{edmonds1970bottleneck}.
	\QEDA

	The result in Theorem~\ref{cluter_blocker} is key to the reformulation of DRBCP-U \eqref{general_eq_dis_robust}, i.e., this result allows us to reformulate max-min-max into max-max-min and to switch the outer maximum operators. It is worthy of noting that for a general set $X$, it can take $O(2^n)$ to find set $\overline{X}$ by enumeration, and it can take up to $O(2^n)$ to find a member of its blocker $F$ according to the definition; however, for special examples in the previous subsection, we have closed-form representations of blocking systems. Thus, we recommend readers to derive closed-form representations if possible.
	
	

	\section{DRBCP-U \eqref{general_eq_dis_robust}: Equivalent Reformulation, MIP Representation, and Confidence Bounds}\label{sec_DRBCP-U}
	
	In this section, we will derive the equivalent reformulation of DRBCP-U \eqref{general_eq_dis_robust}, which motivates us to develop a MIP representation and to study its connection to the sampling average
	approximation counterpart. Similar conclusions follow for DRBCP-U \eqref{general_eq_dis_robust} under $q-$Wasserstein ambiguity set, which can be found in Appendix~\ref{1_norm_resul} except the more computationally expensive reformulations due to the additional dual multiplier.
	
	\subsection{Equivalent Reformulation}
	Using the representation \eqref{eq_dual_wass} of $\infty-$Wasserstein ambiguity set $\P_{\infty}$, we first derive a deterministic representation of DRBCP-U \eqref{general_eq_dis_robust}.
	\begin{theorem}\label{thm_drbcp_d_ref} Suppose $\|\cdot\|=\|\cdot\|_r$ with $r\geq 1$, then DRBCP-U \eqref{general_eq_dis_robust} is equivalent to
		\begin{equation}\label{infity_wasser_re}
		v_U=\frac{1}{N}\sum_{k\in [N]}\max_{y\in F}t_*^k(y),
		\end{equation}
		where $F$ is defined in Theorem \ref{cluter_blocker} and $t_*^k(y)$ is defined as
		\begin{equation}\label{infity_wasser_re_t} 
		t_*^k(y):=\max_{t^k}\left\{ t^k:\sum_{j\in I^k(t^k,y)}(t^k-\bar c_j^k)^r\leq \theta^r\right\}
		\end{equation}
		with set $I^k(t,y)=\{j\in y: t-\overline{c}_j^k\geq 0\}$.
	\end{theorem}
	\proof
	We split the proof into the following two steps.\par
	\noindent\textbf{Step 1. Equivalent Representation of DRBCP-U \eqref{general_eq_dis_robust}.}
	According to the representation \eqref{eq_dual_wass} and Theorem~\ref{cluter_blocker}, we have
	\begin{align*} 
	v_U&=\sup_{\Pr\in \P_{\infty}}\E_{\P}\left[Z(\tilde{\bm{c}})\right] 
	=\frac{1}{N}\sum_{k\in [N]}\max_{\|\bm{c}^k-\bar{\bm{c}}^k\|_r\leq \theta}\max_{y\in F}\min_{j\in y}c_j^k.
	\end{align*}
	By swapping the two maximization operators and linearizing the inner minimization operator, we further have
	\begin{align*} 
	v_U&=\frac{1}{N}\sum_{k\in [N]}\max_{y\in F}\max_{\|\bm{c}^k-\bar{\bm{c}}^k\|_r\leq \theta}\min_{j\in y}c_j^k=\frac{1}{N}\sum_{k\in [N]}\max_{y\in F}\max_{\|\bm{c}^k-\bar{\bm{c}}^k\|_r\leq \theta, t^k\leq c_j^k, \forall j\in y}t^k.
	\end{align*}
	To simplify the inner maximization, we let $\beta_j^k=c_j^k-t^k$ for all $j\in y$ and define set $I^k(t,y)=\{j\in y: t-\overline{c}_j^k\geq 0\}$, then we have 
	\begin{align*}
	v_U
	&=\frac{1}{N}\sum_{k\in [N]}\max_{y\in F}\max_{\bm{\beta}^k,t^k}\left\{t^k:\sum_{j\in y}|t^k+\beta_j^k-\bar{\bm{c}}^k|^r\leq \theta^r, \beta_j^k \geq 0, \forall j\in y\right\}\\
	& =\frac{1}{N}\sum_{k\in [N]}\max_{y\in F}\max_{t^k}\left\{ t^k:\sum_{j\in I^k(t^k,y)}(t^k-\bar c_j^k)^r\leq \theta^r\right\}
	\end{align*}
	where the second equality is due to the optimality condition that at optimality, if $j\in y\setminus I^k(t^k,y) $, we must have $\beta_j^k=-t^k+\bar c_j^k>0$; otherwise, for $j\in I^k(t^k,y) $, given that $t-\overline{c}_j^k\geq 0$ and  $\beta_j^k\geq 0$, the optimal $\beta_j^k=0$.\\
	
	\noindent\textbf{Step 2. Obtaining Optimal $t^k$.} Next, suppose that $y=\left\{y_1,\ldots, y_{|y|}\right\}$ such that $-\infty<\bar{c}_{y_1}^k\leq\ldots\leq \bar{c}_{y_{|\bm{y}|}}^k<\bar{c}_{y_{|\bm{y}|+1}}^k:=+\infty$. Also, suppose that
	\[t_*^k(y)\in \arg\max_{t^k}\left\{ t^k:\sum_{j\in I^k(t^k,y)}(t^k-\bar c_j^k)^r\leq \theta^r\right\}.\]
	Clearly, $t_*^k(y)$ must satisfy
	\begin{subequations}
		\begin{align}\label{infinity_norm_eq_0}
		&\overline{c}^k_{y_{|I^k(t_*^k(y),y)|}}\leq t_*^k(y)<\overline{c}^k_{y_{|I^k(t_*^k(y),y)|+1}},\\\label{infinity_norm_eq_1}
		&\sum_{j\in I^k(t_*^k(y),y)}(\overline{c}^k_{y_{|I^k(t_*^k(y),y)|}}-\bar c_j^k)^r\leq \theta^r,\\\label{infinity_norm_eq_2}
		&\sum_{j\in I^k(t_*^k(y),y)}(\overline{c}^k_{y_{|I^k(t_*^k(y),y)|+1}}-\bar c_j^k)^r> \theta^r.
		\end{align}
	\end{subequations}
	Now, we are going to show that 
	\begin{claim}
		the optimal $t_*^k(y)$ is unique.
	\end{claim}
	\proof To prove the uniqueness of $t_*^k(y)$, it is sufficient to show that the set $I^k(t_*^k(y),y)$ is unique. 
	We prove it by contradiction. Suppose that $I^k(t_*^k(y),y)$ is not unique, then there exists two optimal solutions $t_{1*}^k(y),t_{2*}^k(y)$ such that $I^k(t_{1*}^k(y),y) \subsetneq I^k(t_{2*}^k(y),y) $. Since $t_{1*}^k(y)$ satisfies \eqref{infinity_norm_eq_2}, we have $$\sum_{j\in I^k(t_{1*}^k(y),y)}(\overline{c}^k_{y_{|I^k(t_{1*}^k(y),y)|+1}}-\bar c_j^k)^r> \theta^r.$$
	Since $t_{2*}^k(y)$ satisfies \eqref{infinity_norm_eq_1} and $I^k(t_{1*}^k(y),y) \subsetneq I^k(t_{2*}^k(y),y) $, we have the following two inequalities
	\[\sum_{j\in I^k(t_{2*}^k(y),y)}(\overline{c}^k_{y_{|I^k(t_{2*}^k(y),y)|}}-\bar c_j^k)^r\leq \theta^r,\quad\overline{c}^k_{y_{|I^k(t_{1*}^k(y),y)|+1}}\leq \overline{c}^k_{y_{|I^k(t_{2*}^k(y),y)|}},\]
	which imply that
	\[\theta^r<\sum_{j\in I^k(t_{1*}^k(y),y)}(\overline{c}^k_{y_{|I^k(t_{1*}^k(y),y)|+1}}-\bar c_j^k)^r\leq \sum_{j\in I^k(t_{2*}^k(y),y)}(\overline{c}^k_{y_{|I^k(t_{2*}^k(y),y)|}}-\bar c_j^k)^r\leq \theta^r,\]
	a contradiction.
	\QEDB\par 

	\QEDA
	
	The result in Theorem~\ref{thm_drbcp_d_ref} shows that DRBCP-U \eqref{general_eq_dis_robust} is equivalent to sum of $N$ robust BCPs. Thus, to obtain the value $v_U$ of DRBCP-U \eqref{general_eq_dis_robust}, we can solve $N$ separate robust BCPs, which can be done in parallel. Another observation is that if we choose $L_1$ norm in Theorem~\ref{thm_drbcp_d_ref}, we can further simplify its result.
	\begin{corollary} Suppose $\|\cdot\|=\|\cdot\|_1$ and for each $y\in F$, $y=\left\{y_1,\ldots, y_{|y|}\right\}$ such that $-\infty<\bar{c}_{y_1}^k\leq\ldots\leq \bar{c}_{y_{|\bm{y}|}}^k<\bar{c}_{y_{|\bm{y}|+1}}^k:=+\infty$. Then DRBCP-U \eqref{general_eq_dis_robust} is equivalent to
		\begin{equation}\label{infity_wasser_re_l1}
		v_U=\frac{1}{N}\sum_{k\in [N]}\max_{y\in F}\frac{(\sum_{j\in I^k(y)}\overline{c}_j^k)+\theta}{|I^k(y)|},
		\end{equation}
		where $I^k(y)=\{y_1,\ldots, y_{|I^k(y)|}\}$ and $\overline{c}^k_{y_{|I^k(y)|}}\leq \frac{\sum_{j\in I^k(y)\overline{c}_j^k+\theta}}{|I^k(y)|}<\overline{c}^k_{y_{|I^k(y)|+1}}$. \par 
	\end{corollary}
	According to \eqref{infity_wasser_re_l1}, the value of DRBCP-U is equal to the average of robust values from $N$ scenarios, and within each scenario, we find the subset that maximizes the average of its several smallest elements instead of its bottleneck value. This result will indeed ensure the robustness of the DRBCP-U value.
	
	Finally, we remark the characterization of the worst-case distribution of DRBCP-U \eqref{general_eq_dis_robust}, i.e., the worst-case $\{\bm{c}^k\}_{k\in[N]}$ defined in \eqref{eq_dual_wass}.
	\begin{corollary}\label{cor_wc_pd}For each $k\in [N]$, let $y_*^k\in \arg\max_{y\in F}t_*^k(y)$ and
		\begin{align*}
		c_{j*}^{k}=\begin{cases}
		t_*^k(y_*^k),&\textrm{ if } j\in I^k(t_*^k(y_*^k),y_*^k),\\
		\overline{c}_j^k,&\textrm{ otherwise}.
		\end{cases}
		\end{align*}
		Then the worst-case probability distribution of DRBCP-U \eqref{general_eq_dis_robust} is $\Pr^*=\frac{1}{N}\sum_{k\in [N]}\delta(\tilde{\bm c}-\bm{c}_*^k)$.
	\end{corollary}
	The result in Corollary~\ref{cor_wc_pd} shows that the worst-case distribution has $N$ support points, where each support point slightly deviates from one of the empirical data $\Z=\{\bar{\bm{c}}^k\}_{k\in [N]}$ to enforce the robustness of the DRBCP-U value.
	

	
	\subsection{Mixed-Integer Programming Representation}
	Another side product of Theorem \ref{thm_drbcp_d_ref} is that to compute the optimal value $v_U$ of DRBCP-U \eqref{general_eq_dis_robust} is equivalent to the average of $N$ separate MIPs.
	\begin{proposition}\label{prop_mip_u}
		Suppose that the blocker $F$ admits a binary programming representation $\hat{F}\subseteq \{0,1\}^n$. Then
		$v_U=1/N\sum_{k\in[N]}v_U^k$ and $v_U^k$ is equivalent to
		\begin{align}\label{mixed_integer_form_u}
		v_U^k=\max_{\bm{z}^k\in \hat{F}, t^k, \bm{\beta}^k}&\left\{t^k: \sum_{j\in [n]}|t^k-\bar c_j^k+\beta_j^k|^r\leq \theta^r,\beta_j^k\geq -M_j^k(1-z_j^k), \forall j\in [n]\right\},
		\end{align}
		where $M_{j}^k=\theta+\min_{\tau\in [n]}\bar{c}_\tau^k-\bar{c}_{j}^k$.
	\end{proposition}
	\proof According to the proof of Theorem \ref{thm_drbcp_d_ref}, $v_U=1/N\sum_{k\in[N]}v_U^k$, where
	\[v_U^k:=\max_{y\in F,t^k,\bm{\beta}^k}\left\{t^k:\sum_{j\in y}|t^k-\bar c_j^k+\beta_j^k|^r\leq \theta^r, \beta_j^k\geq 0, \forall j\in y\right\}.\]
	Above, we can augment the vector $\bm{\beta}^k$ by defining free variable $\beta_j^k$ for each $j\in [n]\setminus y$ and rewrite it as
	\[v_U^k:=\max_{y\in F,t^k,\bm{\beta}^k}\left\{t^k:\sum_{j\in [n]}|t^k-\bar c_j^k+\beta_j^k|^r\leq \theta^r, \beta_j^k\geq 0, \forall j\in y, \beta_j^k \in \Re, \forall  j\in [n]\setminus y\right\},\]
	since at optimality, we must have $\beta_{j}^k=-t^k+\bar c_j^k$ for each $j\in [n]\setminus y$.
	
	As the blocker $F$ admits a binary programming representation $\hat{F}$, the value $v_U^k$ is further equal to
	\begin{align}
	v_U^k:=\max_{\bm{z}^k\in \hat{F},t^k,\bm{\beta}^k}\left\{t^k:\sum_{j\in [n]}|t^k-\bar c_j^k+\beta_j^k|^r\leq \theta^r, z_j^k\beta_j^k\geq 0, \forall j\in [n]\right\}.\label{mixed_integer_form_u1}
	\end{align}
	To prove that \eqref{mixed_integer_form_u1} and \eqref{mixed_integer_form_u} are equivalent, we only need to show that $\beta_{j}^k\geq -M_j^k$ for each $j\in[n]$. Since at optimality, $\beta_{j}^k$ is equal to $0$ or $-t^k+\bar c_j^k$, it is sufficient to find an upper bound of $t^k$. 
	According to \eqref{infity_wasser_re_t}, we have
	\[\left(t^k-\min_{\tau\in [n]}\bar{c}_\tau^k\right)^r\leq \sum_{j\in I^k(t^k,y)}(t^k-\bar c_j^k)^r\leq \theta^r,\]
	which implies that $t^k\leq \theta+\min_{\tau\in [n]}\bar{c}_\tau^k$. Thus, we can choose $M_{j}^k=\theta+\min_{\tau\in [n]}\bar{c}_\tau^k-\bar{c}_{j}^k$.
	\QEDA
	
	We remark that: (i) to obtain $v_U$, one needs to solve $N$ MIPs \eqref{mixed_integer_form_u}, which can be done in parallel; and (ii) as long as $r$ is a rational number, using the conic quadratic representation results in \cite{ben2001lectures}, the MIP \eqref{mixed_integer_form_u} essentially can be formulated as a mixed-integer second order conic program (MISOCP).
	
	The following example illustrates the MIP formulation proposed in Proposition~\ref{prop_mip_u}.
	\begin{example}\label{example5} Given a graph $G(V,E)$, suppose that $X$ denotes all $s-t$ paths, then according to Example~\ref{example2}, its blocker $F$ denotes all the $s-t$ cuts. To compute $v_U$ of distributionally robust bottleneck shortest path, it is equal to $1/N\sum_{k\in [N]}v_U^k$, where
		\begin{subequations} \label{norm_eqa_1}
			\begin{align}
			v_U^k=&\max_{\bm{z}^k,\bm{w}^k,\bm{\beta}^k, t^k} \ t^k\\ 
			\text{s.t.} \quad &\|t^k\bm{e}+\bm{\beta}^k-\bar{\bm{c}}^k\|_r\leq \theta, \\
			& \beta_{ij}^k\geq -M_{ij}^k(1-z_{ij}^k), \forall (i,j)\in E, \\
			& z_{ij}^k\leq 1-w_i^k, w_j^k-z_{ij} \leq w_i^k, w_j^k -z_{ij}^k \geq 0, \forall (i,j)\in E,\\ 
			& w_s^k=0,w_t^k=1,w_i^k\in\{0, 1\}, \forall i\in V,
			\end{align}
		\end{subequations}
		where $M_{ij}^k=\theta+\min_{\tau\in E}\bar{c}_\tau^k-\bar{c}_{ij}^k$. 
	\end{example} 
	
	%
	%

	\subsection{Confidence Bounds}
	In this subsection, we compare the optimal value $v_U$ of DRBCP-U \eqref{general_eq_dis_robust} with its sampling average approximation counterpart, which further motivates us on how to choose the Wasserstein radius $\theta$.
	
	Let us define the sampling average approximation of stochastic BCP with empirical distribution $\Pr_{\tilde{\bzeta}}$ as below:
	\begin{subequations}
		\begin{align}
		v_U^{SAA}=\frac{1}{N}\sum_{k\in[N]}\min_{x\in X}\max_{j\in x}\bar{c}_{j}^k.\label{eq_saa}
		\end{align}
		Suppose that the true distribution is $\Pr^{T}$, then the stochastic BCP under the true distribution is defined as
		\begin{align}
		v_U^{T}=\E_{\Pr^T}\left[\min_{j\in X}\max_{j\in h}\tilde{c}_{j}\right].\label{eq_saa_t}
		\end{align}
	\end{subequations}
	Our goal is to determine a proper $\theta$ such that $v_U^{T}+2\theta\geq v_U\geq v_U^{T}$ with a high probability. To begin with, we would like to show the relationship between $v_U$ and $v_U^{SAA}$.
	
	\begin{proposition} \label{infity_norm_saa} Let $v_U^{SAA}$ be defined in \eqref{eq_saa}. Then
		$$\frac{\theta}{{\max}_{y\in F}|y|^{\frac{1}{r}}}\leq v_U-v_U^{SAA}\leq \theta.$$
	\end{proposition}
	\proof 
	\begin{enumerate} [(i)]
		\item To prove $v_U-v_U^{SAA}\leq \theta$, according to Theorem~\ref{cluter_blocker} and Theorem \ref{thm_drbcp_d_ref}, it is sufficient to show that for each $y\in F$, we must have
		\[t_*^k(y)\leq \min_{j\in y} \bar{c}_j^k+\theta.\]
		In \eqref{infity_wasser_re_t}, we have
		\[\left(t_*^k(y)-\min_{j\in y} \bar{c}_j^k\right)^r\leq \sum_{j\in I^k(t_*^k(y),y)}(t_*^k(y)-\bar c_j^k)^r\leq \theta^r,\]
		which implies that $t_*^k(y)\leq \min_{j\in y} \bar{c}_j^k+\theta$.

		\item To prove $v_U-v_U^{SAA}\geq \frac{\theta}{{\max}_{y\in F}|y|^{\frac{1}{r}}}$, we first show that for each $y\in F$, we must have
		\[t_*^k(y)\geq \min_{j\in y} \bar{c}_j^k+\frac{\theta}{|y|^{\frac{1}{r}}}.\]
		In \eqref{infity_wasser_re_t}, we first observe that at optimality, the inequality constraint must be equality:
		\[\sum_{j\in I^k(t_*^k(y),y)}(t_*^k(y)-\bar c_j^k)^r= \theta^r;\]
		otherwise, we can further increase $t_*^k(y)$. Thus,
		\[\theta^r=\sum_{j\in I^k(t_*^k(y),y)}(t_*^k(y)-\bar c_j^k)^r\leq |I^k(t_*^k(y),y)|\left(t_*^k(y)-\min_{j\in y} \bar{c}_j^k\right)^r\leq |y|\left(t_*^k(y)-\min_{j\in y} \bar{c}_j^k\right)^r,\]
		where the first inequality is due to $\bar c_j^k\geq \min_{j\in y} \bar{c}_j^k$ for each $j\in I^k(t_*^k(y),y)$ and the second inequality is because $|I^k(t_*^k(y),y)|\leq |y|$. 
		
		Next, according to Theorem~\ref{cluter_blocker} and Theorem \ref{thm_drbcp_d_ref}, we have
			\begin{align*}
			v_U&=\frac{1}{N}\sum_{k\in [N]}\max_{y\in F}t_*^k(y)\geq \frac{1}{N}\sum_{k\in [N]}\max_{y\in F}\left(\min_{j\in y} \bar{c}_j^k+\frac{\theta}{|y|^{\frac{1}{r}}}\right)\\
			&\geq \frac{1}{N}\sum_{k\in [N]}\max_{y\in F}\min_{j\in y} \bar{c}_j^k+\frac{\theta}{\max_{y\in F}|y|^{\frac{1}{r}}}:= v_U^{SAA}+\frac{\theta}{\max_{y\in F}|y|^{\frac{1}{r}}},
			\end{align*}
			where the last inequality is due to $|y|^{\frac{1}{r}}\leq \max_{h\in F}|h|^{\frac{1}{r}}$ for all $y\in F$.

		This completes the proof.
		\QEDA	
	\end{enumerate}
	
	%
	%
	%
	%
	%

	Now we are ready to show the following main result of the relationship between the value of DRBCP-U $v_U$ and true value $v_U^T$.
	\begin{theorem}\label{thm_bound_u}Suppose that there exists a positive $\sigma$ such that $\E_{\Pr^T}[\exp((Z(\tilde{\bm c})-v_U^T)^2/\sigma^2)]\leq e$ a.s.
		Given $\epsilon\in (0,1)$ and $\theta=N^{-\frac{1}{2}}\sigma \sqrt{-3\log (\epsilon)}{ \max}_{y\in F}|y|^{\frac{1}{r}}=O(N^{-\frac{1}{2}})$, then we have
		\begin{enumerate}[(i)]
			\item 
			$$\Pr^T\left\{v_U\geq v_U^T\right\}\geq 1-\epsilon; \textrm{ and }$$
			\item 
			$$\Pr^T\left\{v_U\leq v_U^T+2\theta\right\}\geq 1-\epsilon.$$
		\end{enumerate}

	\end{theorem}
	\proof Since the proof of Part (ii) is very similar to that of Part (i), we will only prove Part (i). According to Proposition \ref{infity_norm_saa}, it is sufficient to prove
	$$\Pr^T\left\{v_U^{SAA}\geq v_U^T-\frac{\theta}{{ \max}_{y\in F}|y|^{\frac{1}{r}}}\right\}\geq 1-\epsilon.$$
	
	Using Lemma 2 in \cite{guigues2017non} with $d_k=v_U^T-Z(\bar{\bm{c}}^k)$ for each $k\in [N]$, we have
	$$\Pr^T\left\{v_U^{SAA}<v_U^T-\frac{\lambda\sigma}{\sqrt{N}}\right\}\leq e^{-\frac{\lambda^2 }{3}},$$
	i.e.,
	$$\Pr^T\left\{v_U^{SAA}\geq v_U^T-\frac{\lambda\sigma}{\sqrt{N}}\right\}\geq 1-e^{-\frac{\lambda^2 }{3}}.$$
	Letting $e^{-\frac{\lambda^2 }{3}}:=\epsilon$, i.e., $\lambda=\sqrt{-3\log (\epsilon)}$ and $\theta:=\lambda\sigma N^{-\frac{1}{2}}{ \max}_{y\in F}|y|^{\frac{1}{r}}$, we arrive at the conclusion.
	%
	%
	%
	\QEDA 
	
	We make the following remarks for Theorem \ref{thm_bound_u}.
	\begin{enumerate}[(a)]
		\item If we focus on convergence of the optimal objective value instead of probability distribution, then the Wasserstein radius should be chosen in the order of $O(N^{-\frac{1}{2}})$ rather than $O(N^{-\frac{1}{n}})$ shown in \cite{fournier2015rate,esfahani2018data}.
		\item The asymptotic rate $O(N^{-\frac{1}{2}})$ has been observed in \cite{blanchet2019confidence} from a different angle. However, our rate is non-asymptotic.
		\item The presumption of Theorem \ref{thm_bound_u} follows from that of proposition 1 in \cite{guigues2017non}. One can also directly assume that the random vector $\tilde{\bm c}$ is sub-Gaussian. Together with the facts that $|z(\tilde{\bm c})| \leq \max_{j\in [n]}|\tilde{c}_j|$ and maximum of finite sub-Gaussian variables is still sub-Gaussian, the same conclusion follows.
		\item Parts (i) and (ii) together show that with high probability, the value of DRBCP-U $v_U$ is no less than the true value $v_U^T$ and is at most $2\theta$ larger than $v_U^T$. This implies that an appropriate $\theta$ will ensure less conservatism of DRBCP-U value $v_U$.
		\item In the numerical study, we can further improve the radius in Theorem \ref{thm_bound_u} by selecting a list of $\theta$ values and choosing the smallest $\theta$ from this list such that $v_U\geq v_U^T$ with high probability.
	\end{enumerate}

	\section{DRBCP-D \eqref{general_eq_dis_robust_d}: Equivalent Reformulation, and Alternative Decision-Robust Form}\label{sec_DRBCP-D}
	
	In this section, we derive an equivalent reformulation of DRBCP-D \eqref{general_eq_dis_robust_d} and propose an alternative decision-robust model to ensure the robustness of decision.  

	\subsection{Equivalent Reformulation and Confidence Bounds}
	In this subsection, we show that for DRBCP-D \eqref{general_eq_dis_robust_d}, its sampling average approximation counterpart always provides an optimal solution, which motivates us a choice of Wasserstein radius $\theta$.
	
	Let us first define the sampling average approximation of stochastic decision-making BCP as below:
	\begin{subequations}
		\begin{align}
		v_D^{SAA}=\min_{x\in X}\frac{1}{N}\sum_{k\in[N]}\max_{j\in x}\bar{c}_{j}^k.\label{eq_saa_d}
		\end{align}
		Suppose that the true distribution is $\Pr^{T}$, then the stochastic decision-making BCP under the true distribution is defined as
		\begin{align}
		v_D^{T}=\min_{x\in X}\E_{\Pr^T}\left[\max_{j\in x}\tilde{c}_{j}\right].\label{eq_saa_t_d}
		\end{align}
	\end{subequations}
	
	We first observe that DRBCP-D \eqref{general_eq_dis_robust_d} and model \eqref{eq_saa_d} have the same optimal solution.
	\begin{proposition}\label{eq_ref_drbcp_d}
		Let $v_D^{SAA}$ be the optimal value in \eqref{eq_saa_d}. Then
		\[v_D=v_D^{SAA}+\theta.\]
	\end{proposition}
	\proof According to the representation \eqref{eq_dual_wass}, we have
	\begin{align*} 
	v_D&=\min_{x\in X}\sup_{\Pr\in \P_{\infty}}\E_{\Pr}\left[\max_{j\in x}\tilde{c}_{j}\right] 
	=\min_{x\in X}\frac{1}{N}\sum_{k\in [N]}\max_{\|\bm{c}^k-\bar{\bm{c}}^k\|_r\leq \theta}\max_{j\in x}c_j^k=\min_{x\in X}\frac{1}{N}\sum_{k\in [N]}\max_{j\in x}(\bar{c}_j^k+\theta)=v_D^{SAA}+\theta,
	\end{align*}
	where the third equality is due to swapping two maximization operators and the definition of dual norm.
	\QEDA
	
	This result is quite interesting and motivates us to seek an alternative model for DRBCP-D \eqref{general_eq_dis_robust_d}. Besides, we make the following remarks.
	\begin{enumerate}[(i)]
		\item The result in Proposition~\ref{eq_ref_drbcp_d} holds for any $q-$Wasserstein ambiguity set with $q\in [1,\infty)$.
		\item This result shows that applying distributionally robustness into the objective function may indeed only ensure the robustness of objective function. In next subsection, we will propose an alternative decision-robust model for DRBCP-D \eqref{general_eq_dis_robust_d}.
		\item One might want to strengthen the result in Proposition~\ref{eq_ref_drbcp_d} by adding support of the random vector $\tilde{\bm c}$. However, if the support is large enough, we can still obtain the same result.
	\end{enumerate}


	Similar to Proposition \ref{prop_mip_u}, we show that $v_D^{SAA}$ can be obtained by solving an MIP, and according to Proposition~\ref{eq_ref_drbcp_d}, so is $v_D$. The proof is based upon straightforward linearization and is thus omitted.
	\begin{proposition}Suppose that set $X$ admits a binary programming representation $\hat{X}\subseteq \{0,1\}^n$. Then model \eqref{eq_saa_d} is equivalent to the following MIP:
		\begin{equation}\label{formulate_A}
		\begin{aligned}
		\min_{\bm{z}^k\in \hat{X},v_D^k, \forall k\in [N]}\left\{\frac{1}{N} \sum_{k\in[N]}v_D^k:v_D^k \geq (\bar{c}_j^k-L^k)z_j^k+L^k, \forall j\in [n], k\in [N]\right\},
		\end{aligned}
		\end{equation}
		where $L^k:=\min_{\tau\in [n]}\bar{c}_\tau^k$.
	\end{proposition}
	
	
	Next, we will show the relationship between the value of DRBCP-D $v_D$ and the true value $v_D^T$ defined in \eqref{eq_saa_t_d}.
	\begin{theorem}\label{thm_bound_d}Suppose that there exists a positive $\sigma$ such that $\E_{\Pr^T}[\exp((\max_{j\in x}\tilde{c}_{j}-\E_{\Pr^T}[\max_{j\in x}\tilde{c}_{j}])^2/\sigma^2)]\leq e$ for each $x\in X$.
		Given $\epsilon\in (0,1)$, then
		\begin{enumerate}[(i)]
			\item let $\theta=N^{-\frac{1}{2}}\sigma \sqrt{-3\log (\epsilon)+3n\log(2)}=O(N^{-\frac{1}{2}})$. Then we have 
			$$\Pr^T\left\{v_D\geq v_D^T\right\}\geq 1-\epsilon;$$
			\item let $\theta=N^{-\frac{1}{2}}\sigma \sqrt{-3\log (\epsilon)+3n\log(2)}=O(N^{-\frac{1}{2}})$. Then we have 
			$$\Pr^T\left\{v_D\leq v_D^T+2\theta\right\}\geq 1-\epsilon.$$
		\end{enumerate}

	\end{theorem}
	\proof Since the proof of Part (ii) is very similar to that of Part (i), we only prove Part (i). According to Proposition \ref{infity_norm_saa}, it is sufficient to prove
	$$\Pr^T\left\{v_D^{SAA}\geq v_D^T-\theta\right\}\geq 1-\epsilon.$$
	
	Using Lemma 2 in \cite{guigues2017non} with $d_k=\E_{\Pr^T}\left[\max_{j\in x}\tilde{c}_{j}\right]-\max_{j\in x}\bar{c}_{j}^k$ for each $k\in [N]$, we have
	$$\Pr^T\left\{\frac{1}{N}\sum_{k\in[N]}\max_{j\in x}\bar{c}_{j}^k<\E_{\Pr^T}\left[\max_{j\in x}\tilde{c}_{j}\right]-\frac{\lambda\sigma}{\sqrt{N}}\right\}\leq e^{-\frac{\lambda^2 }{3}},$$
	for each $x\in X$. Thus, applying the union bound, we have
	$$\Pr^T\left\{v_D^{SAA}\geq v_D^T-\frac{\lambda\sigma}{\sqrt{N}}\right\}\geq 1-|X|e^{-\frac{\lambda^2 }{3}}\geq 1-2^ne^{-\frac{\lambda^2 }{3}},$$
	where the last inequality is due to $|X|\leq 2^n$. Letting $2^ne^{-\frac{\lambda^2 }{3}}:=\epsilon$, i.e., $\lambda=\sqrt{-3\log (\epsilon)+3n\log(2)}$ and $\theta:=\lambda\sigma N^{-\frac{1}{2}}$, we arrive at the conclusion.
	%
	%
	%
	\QEDA 
	
	In Theorem \ref{thm_bound_d}, we note that {(i) by choosing $\theta=N^{-\frac{1}{2}}\sigma \sqrt{-3\log (\epsilon)}$, we can use the results to estimate the theoretical $(1-2\epsilon)-$ confidence interval of a given solution $x$ as $[1/N\sum_{k\in[N]}\max_{j\in x}\bar{c}_{j}^k-\theta$, $1/N\sum_{k\in[N]}\max_{j\in x}\bar{c}_{j}^k+\theta]$; and (ii) }if the sample size is large enough, then the decision of DRBCP-D \eqref{general_eq_dis_robust_d} will be close to that of the true problem \eqref{eq_saa_t_d}. This might not hold if the sample size is small (i.e., $N$ is small). Thus, to address this issue, in the next subsection, we propose an alternative decision-robust model.
	
	Finally, we remark that the worst-case distribution of DRBCP-D \eqref{general_eq_dis_robust_d}, i.e., the worst-case scenarios $\{\bm{c}^k\}_{k\in[N]}$ are defined in the representation \eqref{eq_dual_wass}.
	\begin{corollary}\label{cor_wc_pd2}Given $x\in X$, for each $k\in [N]$, let $\ell_*^k\in \arg\max_{j\in x}\bar{c}_j^k$ and
		\begin{align*}
		c_{j*}^{k}=\begin{cases}
		\bar{c}_j^k+\theta,&\textrm{ if } j=\ell_*^k,\\
		\bar{c}_j^k,&\textrm{ otherwise}.
		\end{cases}
		\end{align*}
		Then the worst-case probability distribution of DRBCP-D \eqref{general_eq_dis_robust_d} is $\Pr^*=\frac{1}{N}\sum_{k\in [N]}\delta(\tilde{\bm c}-\bm{c}_*^k)$.
	\end{corollary}
	The result in Corollary~\ref{cor_wc_pd2} shows that the worst-case distribution has $N$ support points, where each support point coincides with one of the empirical data $\Z=\{\bar{\bm{c}}^k\}_{k\in [N]}$ except one coordinate.

	

	\subsection{An Alternative Decision-Robust Model}\par 
	%
	\begin{subequations}
		Proposition \ref{eq_ref_drbcp_d} shows that the optimal solution of the model \eqref{eq_saa_d} is always optimal to DRBCP-D \eqref{general_eq_dis_robust_d}, which might not be ideal, especially when the sample size is small. Alternatively, we consider the following lower level set
		\begin{align}
		\L_{\theta}(X)=\left\{x\in X:\frac{1}{N}\sum_{k\in [N]}\max_{j\in x}\bar{c}_j^k\leq v_D:=v_D^{SAA}+\theta \right\},\label{eq_def_L}
		\end{align}
		which corresponds to the set of indifferent solutions to model \eqref{eq_saa_d}. Here, Wasserstein radius $\theta$ provides an indifferent level and needs to be properly chosen (we will discuss this issue at the end of this section). Next, we choose the most robust subset $x$ that minimizes its sampling variance within the lower level set $\L_{\theta}(X)$, i.e., we consider the following ``\textit{decision-robust}'' model:
		\begin{equation} \label{eq_var_for_0_2}
		\begin{aligned}
		v_{DR}=\min_{x\in  \L_{\theta}(X)}\Var_{\Pr_{\tilde{\bzeta}}} \left[\max_{j\in x}\tilde{\zeta}_j\right],
		\end{aligned}
		\end{equation}
		where $\Var_{\Pr} \left[\cdot\right]$ denotes the variance function with probability distribution $\Pr$.
		
	\end{subequations}
	
	Our next result shows that the decision-robust model \eqref{eq_var_for_0_2} can be formulated as an MISOCP.
	\begin{proposition}\label{prop_decision_robust_misocp}Suppose that set $X$ admits a binary programming representation $\hat{X}$. Then the decision-robust model \eqref{eq_var_for_0_2} is equivalent to the following MISOCP
		\begin{equation} \label{var_for}
		\begin{aligned}
		v_{DR}=\min_{\bm{z}\in \hat{X},\bm{v},\bm{\lambda}}\ &\frac{1}{N}\sum_{k\in [N]}\left(v^k-\frac{1}{N}\sum_{\tau\in [N]}v^\tau\right)^2\\
		\text{s.t.}\quad& \frac{1}{N}\sum_{k\in [N]}v^k\leq v_D^{SAA}+\theta,\\
		& v^k \geq (\overline{c}_j^k-L^k)z_j+L^k,\forall \  j \in [n], k \ \in [N],\\
		& v^k \leq  \lambda_j^k\overline{c}_j^k+(1-\lambda_j^k)M^k, \forall \ k \ \in [N],\\
		& \sum_{j\in [n]} \lambda_j^k =1, \forall \  k \in [N],\\
		& \lambda_j^k \leq z_j,  \forall j \ \in [n], \forall \  k \in [N],\\
		& \lambda_j^k \in [0,1],  \forall \  j \in [n], k \ \in [N],
		\end{aligned}
		\end{equation}
		where $L^k:=\min_{\tau\in [n]}\bar{c}_\tau^k$ and $M^k:=\max_{\tau\in [n]}\bar{c}_\tau^k$.
	\end{proposition}
	\proof Since set $X$ admits a binary programming representation $\hat{X}$, according to the definition, the decision-robust model \eqref{eq_var_for_0_2} is equivalent to 
	\begin{equation} \label{eq_var_for_0_3}
	\begin{aligned}
	v_{DR}=\min_{\bm{z}\in \hat{X}}\left\{\frac{1}{N}\sum_{k\in [N]}\left(\max_{j\in [n]}\bar{c}_j^kz_j-\frac{1}{N}\sum_{\tau\in [N]}\left[\max_{j\in [n]}\bar{c}_j^\tau z_j\right]\right)^2:\frac{1}{N}\sum_{k\in [N]}\max_{j\in [n]}\bar{c}_j^kz_j\leq v_D^{SAA}+\theta\right\}.
	\end{aligned}
	\end{equation}
	Next, for each $k\in [N]$, we let $v^k=\max_{j\in [n]}\bar{c}_j^kz_j$ and linearize it by introducing additional variables $\lambda_j^k\in [0,1]$ for all $j\in [n]$. Thus, we arrive at \eqref{var_for}.
	\QEDA
	
	We make the following remarks about decision-robust model \eqref{eq_var_for_0_2}.
	\begin{enumerate}[(i)]
		\item To solve the decision-robust model \eqref{eq_var_for_0_2}, one needs to solve model \eqref{eq_saa_d} first and obtain its optimal value $v_D^{SAA}$.
		\item The decision-robust model \eqref{eq_var_for_0_2} is different from mean-variance optimization \cite{ahmed2006convexity}, where the objective of the latter is to jointly minimize the weighted mean and variance.
		\item The Wasserstein radius $\theta$ here means the indifferent level. One can choose $\theta$ using the results in Theorem \ref{thm_bound_d}. Alternatively, suppose that the optimal solution of model \eqref{eq_saa_d} is $x^{SAA}$, and choose $\theta:=1.645\sqrt{N^{-1}\Var_{\Pr_{\tilde{\bzeta}}} [\max_{j\in x^{SAA}}\tilde{\zeta}_j]}$, which is approximately equal to the upper bound of $95\%$ confidence level of $v_D^{SAA}$. In numerical study, we will use cross-validation method to choose a proper $\theta$.
		\item As long as $\theta\rightarrow 0$ as $N\rightarrow\infty$, the optimal solution to decision-robust model \eqref{eq_var_for_0_2} will eventually converge to the optimal solution of true decision model \eqref{eq_saa_t_d}.
	\end{enumerate}
	
	\begin{subequations}
		The following example shows how to use Proposition~\ref{prop_decision_robust_misocp} to reformulate decision-robust model \eqref{eq_var_for_0_2} as an MISOCP.
		\begin{example}\label{example6} Suppose that there are $m$ ride-sharing service providers and $m$ passengers waiting for a matching, and the travel time $\tilde{c}_{ij}$ from driver $i$ to passenger $j$ is stochastic. Our objective is to find a decision-robust fair matching, which minimizes the largest travel time of the matching. Following Proposition~\ref{prop_decision_robust_misocp}, this decision-robust fair matching can be formulated as
			\begin{equation} \label{bap}\small
			\begin{aligned}
			v_{DR}=\min_{\bm{v, z,\lambda}}\ &\frac{1}{N}\sum_{k\in [N]}\left(v^k-\frac{1}{N}\sum_{\tau\in [N]}v^\tau\right)^2\\
			\text{s.t.}\ & \frac{1}{N}\sum_{k\in [N]}v^k\leq v_D^{SAA}+\theta,\\
			& v^k \geq  ( \overline{c}_{ij}^k-L^k)z_{ij}+L^k, \  \forall i, j\in[m], \forall k \in [N], \\
			& v^k \leq \lambda_{ij}^k \overline{c}_{ij}^k+(1-\lambda_{ij}^k)M^k, \  \forall i, j\in[m], \forall k \in [N], \\
			&  \sum_{i\in [m]}\sum_{j\in [m]} \lambda_{ij}^k =1, \  \forall k \in [N], \\
			&\sum_{i\in [m]}z_{ij}=1, \forall j \in [m],\\
			&\sum_{j\in [m]}z_{ij}=1, \forall i \in [m],\\
			&\lambda_{ij}^k\leq z_{ij},\  \forall i,j\in[m], \forall k \in [N], \\
			& \lambda_{ij}^k \in [0,1], z_{ij}\in\{0,1\} \  \forall i,j\in[m], \forall k \in [N], 
			\end{aligned}
			\end{equation}
			where $L^k=\min_{i,j\in[m]}\bar{c}_{ij}^k,M^k=\max_{i,j\in[m]}\bar{c}_{ij}^k$ and $v_D^{SAA}$ is the optimal value of
			\begin{equation} \label{bap_saa}\small
			\begin{aligned}
			v_D^{SAA}=\min_{\bm{v, z}}\left\{\frac{1}{N}\sum_{k\in [N]}v^k: v^k\geq \overline{c}_{ij}^kz_{ij}, \sum_{i\in [m]}z_{ij}=1, \sum_{j\in [m]}z_{ij}=1,  z_{ij}\in\{0,1\}, \  \forall i, j\in[m],\forall k\in [N]\right\}.
			\end{aligned}
			\end{equation}
		\end{example}
	\end{subequations}

	\section{Numerical Illustrations}\label{sec_numerical}

	In this section, we provide numerical illustrations for both DRBCP-U \eqref{general_eq_dis_robust} and DRBCP-D \eqref{general_eq_dis_robust_d}.  All the instances are coded in Python 3.0 with calls to Gurobi 7.5 on a personal computer with 2.3 GHz Intel Core i5 processor and 8G of memory. 

	
	\subsection{Numerical Illustration of DRBCP-U \eqref{general_eq_dis_robust}}\label{value_infinity}
	
	In this numerical study, we consider a wireless sensor network  (i.e., multi-hop network), where the bandwidth has been normalized to 1 unit. Suppose that there is a channel (Tx-Rx) with transmitter (Tx) and receiver (Rx), which needs to meet a data rate (i.e., minimum capacity) requirement.  Since data rate governs the speed of data transmission of the channel, in order to know how fast we can send data (in bits per second) over the channel, it is important to decide whether this requirement can be guaranteed. 
	We also suppose each data packet is not separable, i.e., there is no flow split, and the link capacity is measured periodically and is broadcast across the network. 
	The uncertainty of link capacity is mainly because of (i) receiver's signal to noise ratio (SNR), and (ii) 
	the network delay. For each channel $i\in [n]$, we use Shannon capacity to calculate the theoretical highest data rate $\tilde{c}_i$ following the work \cite{goldsmith2005wireless}, i.e.,
	\begin{equation*}
	\tilde{c}_i=B_i\log_2\left(1+\textrm{SNR}_i\times \xi_i\right),
	\end{equation*}
	where $B_i$ is the bandwidth, $\xi_i\sim \exp(1)$, $\textrm{SNR}_i$ is the signal-to-noise ratio. Specifically,  $\textrm{SNR}_i$ is computed as	
	\begin{equation*}\textrm{SNR}_i=p_i\sigma^{-2}_i 10^{-12.81-3.76\times \log_{10}{d_{\text{macro}}}},
	\end{equation*}
	where $p_i$ is the transmission power, $\tilde{\sigma}^2$ is the noise level, $d_{i}$ is the distance from the transmitter to the receiver.
	In this numerical study, we suppose that $B_i=1 HZ$, $p_i\sim \textrm{uniform}(0.1, 0.2)W$, $\tilde{\sigma}^2=1 \times 10^{-10} W$, $d_{\text{macro}}\sim \textrm{uniform}(0.03, 0.07)$ km. In addition, we suppose that the network constitutes of 20 nodes and 190 edges, i.e., it is a 20-node complete undirected graph. As the objective of the multi-hop network is to find the smallest data rate, DRBCP-U \eqref{general_eq_dis_robust} becomes $\inf_{\Pr\in \P}\E_{\Pr}\left[\max_{x\in X}\min_{j\in x}\tilde{c}_j\right]$ with $X$ denoting all the $s-t$ paths, where the results in this paper straightforwardly follow.

	To demonstrate the effectiveness of Model \eqref{norm_eqa_1}, we generate three different instances with sample size $N\in \{100,250, 500\}$. For each instance, we consider the Wasserstein radius $\theta\in \{0.00,0.02,\ldots,0.20\}$, norm $\|\cdot\|_1$, and norm $\|\cdot\|_2$ (i.e., $r\in \{1,2\}$). The computational results can be found in Table \ref{bandwidth_num}, where ``Time" represents the total running time for solving an instance with a particular parametric setting, ``$95\%$ of $v_U^{SAA}$" represents the 95\% confidence interval (CI) of the optimal value of Model $\eqref{eq_saa}$, and ``$\theta^*$" is the smallest Wasserstein radius such that its corresponding $v_U$ is below the 95\% asymptotic CI of $v_U^{SAA}$. Here, CI is computed using samples and asymptotic formula, i.e., 
		$$95\% \text{CI}=\text{sample mean $\pm$ 1.96$\times$sample standard deviation/$\sqrt{N}$}.$$
{We also display the theoretical confidence interval, denoted by ``Theoretical 95\% CI", which is computed using the results in Theorem \ref{thm_bound_u} with $\epsilon=0.025$. Note that since we are not able to obtain a closed-form expression of $\sigma$, a key parameter in Theorem \ref{thm_bound_u}, we estimate it by drawing 1000 i.i.d. samples.}


	In Table \ref{bandwidth_num}, we see that by exploring the decomposition structure, each instance can be solved within 3 minutes. On the other hand, there is no much difference between norm $\|\cdot\|_1$ and norm $\|\cdot\|_2$ except the computational time. Thus, in practice, we recommend using $\|\cdot\|_1$ (i.e., $r=1$) in our reformulations. The optimal Wasserstein radius is decreasing as $N$ increases, which is consistent with the intuition that the Wasserstein radius shrinks as the sample size grows. {We also see that the theoretical CI varies as the underlying norm changes, as expected in Theorem \ref{thm_bound_u}. Particularly, both asymptotic CI and theoretical CI are quite similar and do not dominate each other.}


	%
	%
	
	\begin{table}[htbp]
		\begin{center}
			\caption{Numerical Results of DRBCP-U \eqref{general_eq_dis_robust} with Application to the Multi-hop Network. }
			\label{bandwidth_num}
			\small\setlength{\tabcolsep}{4.0pt}
			\begin{tabular}{c|c|ccc|ccc|c|c|c}
				\hline
				\multirow{2}{*}{$N$}  & \multirow{2}{*}{$\theta$} & \multicolumn{3}{c|}{$r=1$}              & \multicolumn{3}{c|}{$r=2$}              & \multirow{2}{*}{$95\%$ CI of $v_U^{SAA}$} &  \multicolumn{2}{c}{{Theoretical $95\%$ CI}} \\ \cline{3-8} \cline{10-11} 
				&                           & Time   & $v_U$ & $\theta^*$             & Time   & $v_U$ & $\theta^*$             &         &{$r=1$ }&{$r=2$ }                                \\ \hline
				\multirow{11}{*}{100} & 0.00                      & 5.210  & 5.411 & \multirow{11}{*}{0.16} & 5.535  & 5.411 & \multirow{11}{*}{0.16} & \multirow{11}{*}{{[}5.279, 5.543{]}} & \multirow{11}{*}{{{[}5.286, 5.554{]}}} &  \multirow{11}{*}{{{[}5.284, 5.369{]}}}    \\
				& 0.02                      & 6.045  & 5.391 &                        & 11.679 & 5.391 &                        &        &  &                                     \\
				& 0.04                      & 6.137  & 5.372 &                        & 12.019 & 5.371 &                        &                                &  &             \\
				& 0.06                      & 6.241  & 5.353 &                        & 12.999 & 5.352 &                        &                          &  &                   \\
				& 0.08                      & 6.222  & 5.335 &                        & 12.970 & 5.332 &                        &                                  &  &           \\
				& 0.10                      & 6.179  & 5.316 &                        & 12.318 & 5.313 &                        &                                     &  &        \\
				& 0.12                      & 6.260  & 5.299 &                        & 12.714 & 5.294 &                        &                             &  &                \\
				& 0.14                      & 6.290  & 5.281 &                        & 12.889 & 5.275 &                        &                             &  &                \\
				& 0.16                      & 6.290  & 5.264 &                        & 12.317 & 5.257 &                        &                                  &  &           \\
				& 0.18                      & 6.285  & 5.247 &                        & 12.658 & 5.238 &                        &                                    &  &         \\
				& 0.20                      & 6.320  & 5.230 &                        & 12.355 & 5.220 &                        &                                    &  &         \\ \hhline{===========}
				\multirow{11}{*}{250} & 0.00                      & 14.057 & 5.417 & \multirow{11}{*}{0.10} & 14.292 & 5.417 & \multirow{11}{*}{0.10} & \multirow{11}{*}{{[}5.335, 5.500{]}}  & \multirow{11}{*}{{{[}5.335, 5.504{]}}} &\multirow{11}{*}{{{[}5.337, 5.391{]}}}    \\
				& 0.02                      & 16.462 & 5.398 &                        & 31.270 & 5.397 &                        &                                 &  &            \\
				& 0.04                      & 16.473 & 5.378 &                        & 31.542 & 5.378 &                        &                                     &  &        \\
				& 0.06                      & 16.762 & 5.359 &                        & 33.610 & 5.358 &                        &                                     &  &        \\
				& 0.08                      & 16.823 & 5.341 &                        & 33.901 & 5.339 &                        &                                    &  &         \\
				& 0.10                      & 16.734 & 5.322 &                        & 32.767 & 5.319 &                        &                                   &  &          \\
				& 0.12                      & 16.978 & 5.304 &                        & 33.183 & 5.300 &                        &                                   &  &          \\
				& 0.14                      & 16.939 & 5.286 &                        & 33.352 & 5.281 &                        &                                   &  &          \\
				& 0.16                      & 16.999 & 5.268 &                        & 33.135 & 5.262 &                        &                                  &  &           \\
				& 0.18                      & 16.982 & 5.251 &                        & 33.149 & 5.244 &                        &                                   &  &          \\
				& 0.20                & 17.099 & 5.234 &                              & 32.443 & 5.225 &                        &                                  &  &         \\ \hhline{===========}
				\multirow{11}{*}{500} & 0.00                      & 29.429 & 5.452 & \multirow{11}{*}{0.06} & 32.524 & 5.452 & \multirow{11}{*}{0.06} & \multirow{11}{*}{{[}5.397, 5.508{]}}  & \multirow{11}{*}{{{[}5.393, 5.513{]}}}& \multirow{11}{*}{{{[}5.395, 5.433{]}}}     \\
				& 0.02                      & 34.224 & 5.433 &                        & 65.434 & 5.432 &                        &                               &  &              \\
				& 0.04                      & 34.619 & 5.413 &                        & 66.014 & 5.413 &                        &                                   &  &          \\
				& 0.06                      & 35.205 & 5.394 &                        & 68.920 & 5.393 &                        &                                &  &             \\
				& 0.08                      & 35.266 & 5.375 &                        & 69.517 & 5.374 &                        &                                &  &             \\
				& 0.10                      & 34.839 & 5.357 &                        & 67.238 & 5.354 &                        &                                &  &             \\
				& 0.12                      & 35.302 & 5.339 &                        & 68.342 & 5.335 &                        &                                &  &             \\
				& 0.14                      & 35.585 & 5.321 &                        & 69.470 & 5.316 &                        &                                &  &             \\
				& 0.16                      & 35.526 & 5.303 &                        & 70.626 & 5.297 &                        &                                &  &             \\
				& 0.18                      & 35.470 & 5.286 &                        & 69.109 & 5.279 &                        &                                 &  &            \\
				& 0.20                      & 35.668 & 5.269 &                        & 68.809 & 5.260 &                        &                                 &  &            \\ \hline
			\end{tabular}
			\vspace*{-10pt}
		\end{center}
	\end{table}



	

	Finally, we also test DRBCP-U \eqref{general_eq_dis_robust} under $q-$Wasserstein ambiguity set with $q=2$ using the same instances, where the formulation can be found in the appendix. We display the numerical results in Table \ref{dro_value_w1} in Appendix~\ref{appendix_num_2w} due to the page limit, where we let $\lambda^*$ denote the best dual multiplier found. For norm $\|\cdot\|_2$ (i.e., $r=2$), each instance takes more than an hour to solve, and is thus not included in Table \ref{dro_value_w1}. Compared to the results in Table~\ref{bandwidth_num} , it is seen that the running time of DRBCP-U \eqref{general_eq_dis_robust} under $\infty-$Wasserstein ambiguity set is much shorter than that of DRBCP-U \eqref{general_eq_dis_robust} under $2-$Wasserstein ambiguity set, which is mainly because the reformulation of the former does not have a dual multiplier. In contrast, the reformulation of the latter has. On the other hand, we can also see that by choosing a proper Wasserstein radius, DRBCP-U under each ambiguity set delivers a similar objective value, respectively. Thus, this implies that in practice, we recommend using the $\infty-$Wasserstein ambiguity set due to its shorter running time and similar statistical performances.

	\subsection{Numerical Illustrations of Decision Robust Model \eqref{eq_var_for_0_2}}\label{sec_num_decision}
	
	We generate two sets of numerical study to illustrate the decision robust model \eqref{eq_var_for_0_2}.
	\vspace{5pt}
	
\begin{figure}[htbp]
		\centering
		\subfloat[0.48\textwidth][Mean v.s. $\theta$]{
			\includegraphics[width=0.48\textwidth]{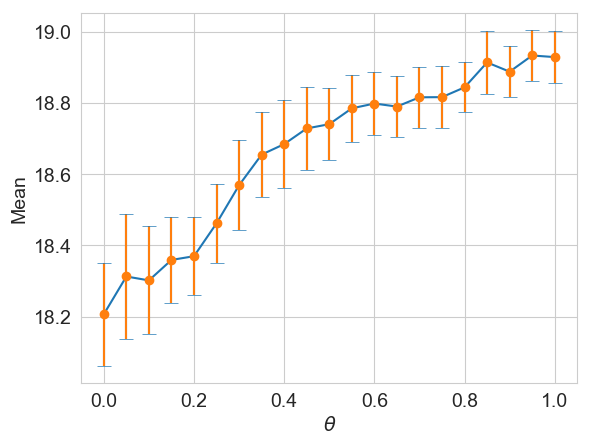} 
			\label{mean_50}
		}\hfill
		\subfloat[0.48\textwidth][Variance v.s. $\theta$]{
			\centering
			\includegraphics[width=0.48\textwidth]{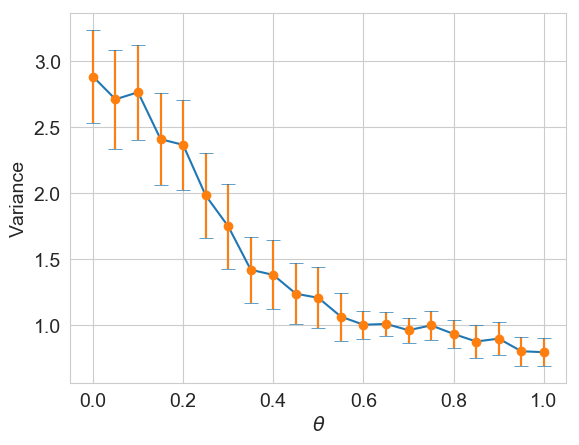} 
			\label{var_50}}
		\caption{Numerical Illustration of Decision Robust Model \eqref{eq_var_for_0_2}: Using Real-world Dataset}\label{fig_num_2}
	\end{figure}

	\noindent{\textbf{First Set of Numerical Study: Using Real-world Dataset-}} To demonstrate the effectiveness of decision robust model \eqref{eq_var_for_0_2}, in our first set of numerical study, we test the bottleneck matching Example \ref{example6}. We suppose that there are 9 drivers and 9 passengers, where the travel time between each pair of driver and passenger is stochastic. We use the 2018 New York taxi record data from \url{https://www1.nyc.gov/site/tlc/about/tlc-trip-record-data.page}, where we treat average monthly travel time as a sample and thus, there are 12 empirical samples in total. 
	We use the following cross-validation method to choose the optimal $\theta$: (i) We randomly choose 8 samples as the training dataset and the other 4 samples as the testing dataset; (ii) Given each $\theta\in \{0, 0.05, 0.1, \cdots, 1.0\}$, we solve the decision robust model \eqref{bap} using the training dataset and evaluate the mean and variance of its unique optimal solution using the testing dataset; and (iii) We repeat steps (i) and (ii) 50 times and compute the $95\%$ confidence intervals of the estimated mean and variance. The computational results are illustrated in Figure~\ref{fig_num_2}. 
	
	From Figure~\ref{fig_num_2}, it is seen that the expected objective value slightly increases as $\theta$ grows. On the other hand, the variance reduces significantly as $\theta$ increases, in particular, when $\theta$ is small. Thus, in general, we recommend users choosing a proper $\theta$ such that the objective value of its corresponding optimal solution is relatively the same as SAA, but its variability is much smaller. In this numerical instance, we recommend using $\theta=0.65$, where the expected objective value increases only by $3.3\%$, but its variance reduces by $64.2\%$. 
	
We also study distributionally robust fair matching problem in Example \ref{example6} under total variation based ambiguity set (see theorem 1 in \cite{jiang2018risk}), the detailed formulation can be found in Appendix \ref{appendix_num_phi} and the numerical results are shown in Figure~\ref{fig_num_2_div}. As observed in \cite{lam2019recovering,duchi2019variance}, the formulation using total variation (e.g., a special case of $\phi$-divergence) based ambiguity set (e.g., a special case of $\phi$-divergence) can help reduce the variance of the decision. However, comparing Figure~\ref{fig_num_2} with Figure~\ref{fig_num_2_div},  it is seen that our decision robust model tends to find solutions with much lower variances while the objective value is almost the same as that of SAA.
	
	\vspace{5pt}
	
	\noindent{\textbf{Second Set of Numerical Study: Using Hypothetical Dataset}} In our second set of numerical study, 
	similar to the first set of numerical study, we consider 9 drivers and 9 passengers. Their travel times are independent Gaussian variables truncated at the interval $[0,\infty)$, where their means are equal to the average travel times of the 12 empirical samples of 2018 New York record data, denoted by $\bm\mu_{NY}$, and their standard deviations are proportional to the empirical standard deviations of the 2018 New York data, $\alpha\bm\sigma_{NY}$ with $\alpha>1$. In particular, we consider two cases with $\alpha=10$ and $\alpha=100$. For each case, we generate three different testing instances with sample size $N\in\{10, 20, 40\}$, respectively. For the case $\alpha=10$, we consider the candidate $\theta$ ranging from $\{0, 1, \cdots, 10\}$, while for the case $\alpha=100$, we consider the candidate $\theta$ ranging from $\{0, 5, \cdots, 50\}$. 
	In addition, for each instance, we evaluate the unique optimal solution of the decision robust model \eqref{bap} and its 95\% confidence interval by generating $10^5$ i.i.d. samples. {We also compute the theoretical confidence interval, denoted by ``Theoretical 95\% CI", according to Theorem \ref{thm_bound_d} with $\epsilon=0.025$. Note that as  we are not able to obtain a closed-form expression of $\sigma$, a key parameter in Theorem \ref{thm_bound_d}, we estimate it using the same $10^5$ i.i.d. samples.} 
 The computational results are displayed in Table \ref{uber_match}.

	\begin{table}[htbp]
		\begin{center}
			\caption{Numerical Results of Decision Robust Model \eqref{eq_var_for_0_2}: Using Hypothetical Dataset }
			\label{uber_match}
			\scriptsize\setlength{\tabcolsep}{4.0pt}
			{ 
			\begin{tabular}{c|c|c|c|c|c|c|c|c|c|c}
				\hline
				\multirow{2}{*}{$N$}    & \multicolumn{5}{c|}{Case 1: $\N\left(\bm\mu_{NY},  10^2\bm\sigma_{NY}\right)$}                                                   & \multicolumn{5}{c}{Case 2: $\N\left(\bm\mu_{NY},  10^4\bm\sigma_{NY}\right)$}                                                       \\ \cline{2-11} 
				& $\theta$ & Time     & $v_{DR}$                    & $95\%$ CI  of  Solutions  &Theoretical $95\%$ CI       & $\theta$ & Time    & $v_{DR}$                      & $95\%$ CI  of  Solutions    &Theoretical $95\%$ CI            \\ \hline
				\multirow{11}{*}{10} & 0        & 2.181    & 345.082                  & $66.646\pm0.130$      & $66.646\pm 0.220
				$              & 0       & 2.027   & \multirow{2}{*}{25974.710} & \multirow{2}{*}{$497.251\pm 1.382$}& \multirow{2}{*}{$497.251\pm 2.345$} \\ \cline{2-8}
				& 1        & 2.335    & 266.640                  & $66.198\pm0.139$           & $66.198 \pm 0.235$              & 5       & 1.923   &                            &                                        \\ \cline{2-11} 
				& 2        & 1.768    & \multirow{2}{*}{128.644} & \multirow{2}{*}{$68.351\pm 0.127$}  & \multirow{2}{*}{$68.351\pm0.216$} 
				 & 10      & 1.750   & \multirow{2}{*}{19726.124} & \multirow{2}{*}{$478.010\pm 1.253$} & \multirow{2}{*}{$478.010\pm 2.127$}    
				  \\ \cline{2-3} \cline{7-8}
				& 3        & 2.076    &                          &                        &            & 15      & 1.419   &                            &                                        \\ \cline{2-11} 
				& 4        & 1.805    & 86.938                   & $63.478\pm 0.130$    & $63.478\pm 0.220$                  & 20      & 1.475   & \multirow{3}{*}{15999.190} & \multirow{3}{*}{$489.546\pm  1.236$}  & \multirow{3}{*}{$489.546\pm  2.097$}\\ \cline{2-8}
				& 5        & 1.518    & \multirow{3}{*}{84.734}  & \multirow{3}{*}{$66.635\pm 0.127$} & \multirow{3}{*}{$66.635\pm 0.215$} & 25      & 1.635   &                            &                                        \\ \cline{2-3} \cline{7-8}
				& 6        & 1.310    &                          &             &                       & 30      & 2.121   &                            &                                        \\ \cline{2-3} \cline{7-11} 
				& 7        & 1.416    &                          &                &                    & 35      & 1.851   & 12968.835                  & $473.883\pm  1.242$    & $473.883\pm  2.108$                \\ \cline{2-11} 
				& 8        & 1.698    & 76.076                   & $63.968\pm 0.128$     & $63.968\pm 0.217$                  & 40      & 1.656   & 12879.386                  & {\bm{$453.492\pm  1.251$}}      & {{$453.492\pm 2.123$}}                \\ \cline{2-11} 
				& 9        & 1.458    & 66.968                   & {\bm{$60.016\pm 0.119$} }   & {{$60.016\pm 0.203$} }               & 45      & 1.267   & 12303.694                  & $463.123\pm 1.246$   & $463.123\pm2.115$                    \\ \cline{2-11} 
				& 10       & 1.598    & 52.736                   & $67.945\pm 0.131$     & $67.945\pm 0.222$                & 50      & 1.447   & 12088.447                  & $474.230\pm  1.226$      & $474.230\pm  2.081$                 \\ \hhline{===========}
				\multirow{11}{*}{20} & 0        & 5.975    & \multirow{2}{*}{129.040} & \multirow{2}{*}{$64.258\pm 0.165$}& \multirow{2}{*}{$64.258\pm 0.280$} 
				& 0       & 6.483   & \multirow{3}{*}{17812.189} & \multirow{3}{*}{$427.099\pm  1.149$}&\multirow{3}{*}{$427.099\pm 1.950$} \\ \cline{2-3} \cline{7-8}
				& 1        & 5.767    &                          &               &                     & 5       & 8.392   &                            &                                        \\ \cline{2-8}
				& 2        & 4.982    & \multirow{6}{*}{108.433} & \multirow{6}{*}{$61.118\pm 0.159$}& \multirow{6}{*}{$61.118\pm 0.271$} & 10      & 6.694   &                            &                                        \\ \cline{2-3} \cline{7-11} 
				& 3        & 43.247   &                          &                      &          & 15      & 54.704  & \multirow{2}{*}{14899.214} & \multirow{2}{*}{$425.33\pm  1.159$}  & \multirow{2}{*}{$425.33\pm  1.968$}    \\ \cline{2-3} \cline{7-8}
				& 4        & 5.733    &                          &                      &              & 20      & 7.281   &                            &                                        \\ \cline{2-3} \cline{7-11} 
				& 5        & 84.842   &                          &                   &                 & 25      & 7.193   & 12579.868                  & $406.824\pm 1.030$  & $406.824\pm 1.749$                    \\ \cline{2-3} \cline{7-11} 
				& 6        & 7.103    &                          &                     &               & 30      & 45.405  & \multirow{2}{*}{12433.038} & \multirow{2}{*}{$473.365\pm 1.644$}  & \multirow{2}{*}{$473.365\pm  2.790$}  \\ \cline{2-3} \cline{7-8}
				& 7        & 7.136    &                          &               &                     & 35      & 34.292  &                            &                                        \\ \cline{2-11} 
				& 8        & 63.746   & 106.171                  & $63.520\pm 0.169$     & $63.520\pm 0.287$                & 40      & 6.049   & \multirow{3}{*}{12268.799} & \multirow{3}{*}{{\bm{$397.311\pm 1.007$} }} & \multirow{3}{*}{{{$397.311\pm 1.710$} }}  \\ \cline{2-8}
				& 9        & 9.318    & 102.174                  & $63.133\pm 0.120$    & $63.133\pm 0.203$                   & 45      & 59.523  &                            &                                        \\ \cline{2-8}
				& 10       & 9.459    & 89.397                   & {\bm{$60.904\pm 0.122$}}    & {{$60.904\pm 0.207$}}                  & 50      & 50.737  &                            &                                        \\\hhline{===========}
				\multirow{11}{*}{40} & 0        & 83.605   & \multirow{3}{*}{166.063} & \multirow{3}{*}{$61.118\pm 0.159$} & \multirow{3}{*}{$61.118\pm 0.271$} & 0       & 146.536 & 29951.644                  & $412.923\pm  1.086$            & $412.923\pm   1.843$           \\ \cline{2-3} \cline{7-11} 
				& 1        & 66.774   &                          &                      &              & 5       & 123.786 & 21742.854                  & $412.998\pm 1.063$    & $412.998\pm 1.805$                  \\ \cline{2-3} \cline{7-11} 
				& 2        & 94.462   &                          &                      &              & 10      & 139.989 & \multirow{7}{*}{16679.337} & \multirow{7}{*}{$397.311\pm  1.007$}  & \multirow{7}{*}{$397.311\pm 1.710$}  \\ \cline{2-8}
				& 3        & 110.405  & \multirow{4}{*}{148.079} & \multirow{4}{*}{{\bm{$56.372\pm 0.103$}}}  & \multirow{4}{*}{{{$56.372\pm 0.174$}}} & 15      & 118.074 &                            &                                        \\ \cline{2-3} \cline{7-8}
				& 4        & 33.100   &                          &                           &         & 20      & 151.350 &                            &                                        \\ \cline{2-3} \cline{7-8}
				& 5        & 160.091  &                          &                            &       & 25      & 130.666 &                            &                                        \\ \cline{2-3} \cline{7-8}
				& 6        & 282.649  &                          &                     &               & 30      & 141.449 &                            &                                        \\ \cline{2-8}
				& 7        & 293.787  & \multirow{3}{*}{133.785} & \multirow{3}{*}{$63.109\pm 0.120$}& \multirow{3}{*}{$63.109\pm 0.204$} & 35      & 135.817 &                            &                                        \\ \cline{2-3} \cline{7-8}
				& 8        & 511.723  &                          &                       &             & 40      & 136.939 &                            &                                        \\ \cline{2-3} \cline{7-11} 
				& 9        & 1263.980 &                          &                       &             & 45      & 193.149 & \multirow{2}{*}{13742.891} & \multirow{2}{*}{{\bm{$396.782\pm  1.010$}}}& \multirow{2}{*}{{$396.782\pm  1.715$}}  \\ \cline{2-8}
				& 10       & 49.186   & 148.760                  & $67.394\pm 0.125$       & $67.394\pm 0.244$               & 50      & 160.060 &                            &                                        \\ \hline
			\end{tabular}
		}
		\end{center}
	\end{table}

	In Table \ref{uber_match}, ``Time" represents the total running time for solving an instance, ``$95\%$ CI of Solutions" represents the 95\% confidence interval (CI) of evaluating the optimal solution of the decision robust model \eqref{bap}. Note that $\theta=0$ represents SAA model. It is seen that the SAA solution can be inferior to the solution of the decision robust model \eqref{bap}, even for the instances with the larger sample size $N=40$. The phenomenon is even more striking when the variance of the true distribution increases. Since the deviations after positive-negative sign in the column ``$95\%$ CI of Solutions" are proportional to the sample standard deviation of $10^5$ samples, it is seen that for most of the instances with $\theta>0$, their confidence intervals are smaller than SAA model, which is consistent with the objective of the decision robust model \eqref{bap}. On the other hand, since the decision robust model \eqref{bap} involves big-M coefficients and more auxiliary variables, its running time is often longer than that of the SAA. {It is seen that the theoretical CI is often wider than the asymptotic one, implying that the proposed theoretical results might not be tight.} In practice, if the sample size is small and the distributionally robust optimization model does not offer better solutions, we recommend using decision robust model as a better alternative.  
	
	We also test distributionally robust fair matching problem in Example \ref{example6} under total variation based ambiguity set in Appendix \ref{appendix_num_phi} and the numerical results are shown in Table~\ref{uber_match_div}. Comparing Table~\ref{uber_match} with Table~\ref{uber_match_div}, we can see that for most of the instances, our decision robust model can find better decisions with lower cross-validation objective values than those of distributionally robust counterpart under total variation based ambiguity set. This further demonstrates that our proposed decision-robust model proposed in Section~\ref{sec_DRBCP-D} can be indeed more reliable than the alternatives.

	\section{Extension: Distributionally Robust $\Gamma$-sum Bottleneck Combinatorial Problems (DR$\Gamma$BCP)}\label{sec_T_sum_dro}
	In this section, we extend the developed results into distributionally robust $\Gamma$-sum bottleneck combinatorial problems (DR$\Gamma$BCP), where the decision-makers are interested in the worst-case sum of $\Gamma-$ largest elements instead of the bottleneck element. Thus, following the notation of BCP \eqref{eq_cbp}, DRBCP-U \eqref{general_eq_dis_robust}, and DRBCP-D \eqref{general_eq_dis_robust}, their  $\Gamma$-sum counterparts, i.e., $\Gamma$BCP, DR$\Gamma$BCP-U, DR$\Gamma$BCP-U, can be formulated as below:
	\begin{subequations}
		\begin{align}
		\textrm{($\Gamma$BCP)}\quad &Z_{\Gamma}(\bm{c})=\min_{x\in X}\max_{s\subseteq x, |s|=\Gamma}\sum_{j\in s}c_j, \label{eq_cbp_gamma}\\
		\textrm{(DR$\Gamma$BCP-U)}\quad&v_{U\Gamma }=\sup_{\Pr\in \P}\E_{\Pr}\left[Z_{\Gamma}(\tilde{\bm{c}})\right],\label{general_eq_dis_robust_gamma}\\
		\textrm{(DR$\Gamma$BCP-D)}\quad &v_{D\Gamma }=\min_{x\in X}\sup_{\Pr\in \P}\E_{\Pr}\left[\max_{s\subseteq x, |s|=\Gamma}\sum_{j\in s}\tilde c_j\right],\label{general_eq_dis_robust_gamma_d}
		\end{align}
	\end{subequations}
	where we assume that all the elements in $X$ have the size at least $\Gamma$. Please note that in \cite{natarajan2013probabilistic}, the authors studied DR$\Gamma$BCP-D with known marginal distributions, while this paper studies both DR$\Gamma$BCP-U and DR$\Gamma$BCP-D under the less conservative Wasserstein ambiguity set.
	
	\subsection{DR$\Gamma$BCP-U \eqref{general_eq_dis_robust_gamma}: Equivalent Reformulations and Confidence Bounds}
	
	Following the same definition of $\bar{X}$ in Section \ref{sec_dual_bcp}, we can show that the following result. The proof is similar to Lemma~\ref{cluter_blocker_0} and is thus omitted.
	\begin{lemma}\label{cluter_blocker_0_gamma}$\Gamma$BCP \eqref{eq_cbp_gamma} is equivalent to
		\begin{equation} \label{equa_clust_gamma}
		Z_\Gamma(\bm{c})=\min_{h\in \overline{X}}\max_{s\subseteq h, |s|=\Gamma}\sum_{j\in s}c_j.
		\end{equation}
	\end{lemma}

	Next, let us define $\Gamma-$blocker $F_{\Gamma}$ as
	\[F_{\Gamma}=\left\{y\subseteq 2^{[N]}: \forall h\in\overline{X} \ \text{ such that $y\cap {h\choose \Gamma}\neq \emptyset$ and}\ y \ \text{is minimal with this property}\right\},\]
	where for any given set $S$, we let ${S\choose\Gamma}$ denote the collection of all its size-$\Gamma$ subsets. According to its definition, there exists a unique $\Gamma-$blocker $F_{\Gamma}$, where the existence is because ${[N]\choose \Gamma}$ satisfies the definition, and the uniqueness is due to the minimality property. Another observation is that we must have $F_T\subseteq {[N]\choose \Gamma}$ according to the nonempty intersection requirement and minimality property.
	Then we present the equivalent dual formulation of $\Gamma$BCP \eqref{equa_clust_gamma}.
	\begin{theorem}\label{cluter_blocker_gamma}Suppose $F$ is the $\Gamma-$blocker of $\overline{X}$. Then $\Gamma$BCP \eqref{eq_cbp_gamma} is equivalent to
		\begin{equation}\label{block_gamma}
		Z_\Gamma(\bm{c})=\max_{y\in F_{\Gamma}}\min_{s\in y}\sum_{j\in s}c_j.
		\end{equation}
	\end{theorem}
	\proof Let $\hat{Z}_\Gamma(\bm{c})$ denote the value of the right-hand side of \eqref{block_gamma}. We separate the proof into two steps.
	\begin{enumerate}[(i)]
		\item $Z_\Gamma(\bm{c})\geq \hat{Z}_\Gamma(\bm{c})$. For any $y\in F_{\Gamma}$ and $x\in \bar{X}$, we have ${x\choose \Gamma}\cap y\neq \emptyset$. Thus,
		\begin{align*}
		\min_{s\in y}\sum_{j\in s}c_j\leq \min_{s\in {x\choose \Gamma}\cap y}\sum_{j\in s}c_j \leq \max_{s\in {x\choose \Gamma}\cap y}\sum_{j\in s}c_j \leq \max_{s\in {x\choose \Gamma}}\sum_{j\in s}c_j,
		\end{align*}
		where the first inequality is due to ${x\choose \Gamma}\cap y\subseteq y$, the second inequality is due to the fact that minimization is no larger than maximization, and the third one is due to ${x\choose \Gamma}\cap y\subseteq {x\choose \Gamma}$. Since the above inequalities hold for any $y\in F_{\Gamma}$ and $x\in \bar{X}$, we must have $Z_\Gamma(\bm{c})\geq \hat{Z}_\Gamma(\bm{c})$.
		
		\item $Z_\Gamma(\bm{c})\leq \hat{Z}_\Gamma(\bm{c})$. Let set $\{s_i\}_{i\in [\binom{n}{\Gamma}]}:={[n]\choose \Gamma}$ be such that their corresponding objective values are sorted in ascending order, i.e., $\sum_{j\in s_1}c_j\leq \ldots\leq \sum_{j\in s_i}c_j\leq \ldots \leq \sum_{j\in s_{\binom{n}{\Gamma}}}c_j$. 
		Let $i_1\in [\binom{n}{\Gamma}]$ be the smallest index such that ${x\choose \Gamma}\subseteq \{s_i\}_{i\in [i_1]}$ for some $x\in \bar{X}$. Clearly, for any such $x$, according to our construction, we must have $s_{i_1}\subseteq x$, 
		otherwise, we can reduce $i_1$, a contradiction. Thus, we have
		\[Z_\Gamma(\bm{c})=\sum_{j\in s_{i_1}}c_j.\]
		
		On the other hand, since ${x\choose \Gamma}\not\subseteq \{s_i\}_{i\in [i_1-1]}$ for any $x\in \bar{X}$, for any $x\in \bar{X}$, there must exist $s_{ix}\subseteq x$ for some $ix\in \{i_1,\ldots, \binom{n}{\Gamma}\}$. Therefore, according to the minimality property of the elements of $\Gamma-$blocker $F_{\Gamma}$, there exists at least one $\hat{y}\in F_{\Gamma}$ such that $\hat{y}\subseteq \{s_i\}_{i\in [i=i_1,\ldots, \binom{n}{\Gamma}]}$. Hence,
		\[\hat{Z}_\Gamma(\bm{c})=\max_{y\in F_{\Gamma}}\min_{s\in y}\sum_{j\in s}c_j \geq \min_{s\in \hat y}\sum_{j\in s}c_j \geq\sum_{j\in s_{i_1}}c_j=Z_\Gamma(\bm{c}).\]
		\QEDA
	\end{enumerate}
	
	We are ready to present the main results of reformulations and confidence bounds. 
	\begin{subequations}
		\begin{theorem}[Reformulations]\label{thm_drbcp_d_ref_gamma} Suppose $\|\cdot\|=\|\cdot\|_r$ with $r\geq 1$, then
			\begin{enumerate}[(i)]
				\item DR$\Gamma$BCP-U \eqref{general_eq_dis_robust_gamma} is equivalent to
				\begin{equation}\label{infity_wasser_re_gamma}
				v_{U\Gamma}=\frac{1}{N}\sum_{k\in [N]}\max_{y\in F_\Gamma}\max_{\bm{\beta}^k}\left\{\min_{s\in y}\sum_{j\in s}(\bar{c}_j^k+\beta_j^k):\sum_{j\in \bigcup_{s\in y}s }(\beta_j^k)^r\leq \theta^r, \beta_j^k\geq 0,\forall j\in \bigcup_{s\in y}s \right\}.
				\end{equation}
				
				\item Let $v_{U\Gamma}^{SAA}=\frac{1}{N}\sum_{k\in[N]}\min_{x\in X}\max_{s\subseteq x, |s|=\Gamma}\sum_{j\in s}\bar{c}_{j}^k$ denote the sampling average approximation counterpart of DR$\Gamma$BCP-U \eqref{general_eq_dis_robust_gamma}. Then
				\begin{equation}\label{infity_wasser_re_gamma_bound}
				v_{U\Gamma}^{SAA}+\frac{\Gamma\theta}{\max_{y\in F_{\Gamma}}|\bigcup_{s\in y}s |^{\frac{1}{r}}}\leq v_{U\Gamma}\leq v_{U\Gamma}^{SAA}+\Gamma^{\frac{r-1}{r}}\theta.
				\end{equation}
			\end{enumerate} 
	\end{theorem}\end{subequations}
	\proof
	\begin{enumerate}[(i)]
		\item 
		According to the representation \eqref{eq_dual_wass} and Theorem~\ref{cluter_blocker_gamma}, we have
		\begin{align*} 
		v_{U\Gamma}&=\sup_{\Pr\in \P_{\infty}}\E_{\P}\left[Z_\Gamma(\tilde{\bm{c}})\right] 
		=\frac{1}{N}\sum_{k\in [N]}\max_{\|\bm{c}^k-\bar{\bm{c}}^k\|_r\leq \theta}\max_{y\in F_\Gamma}\min_{s\in y}\sum_{j\in s}c_j^k.
		\end{align*}
		By swapping the two maximization operators and letting $\beta_j^k=c_j^k-\bar{c}_j^k$, we further have
		\begin{align*} 
		v_{U\Gamma}&=\frac{1}{N}\sum_{k\in [N]}\max_{y\in F_\Gamma}\max_{\|\bm{\beta}^k\|_r\leq \theta}\min_{s\in y}\sum_{j\in s}(\bar{c}_j^k+\beta_j^k).
		\end{align*}
		According to the optimality condition, we must have $\beta_j^k\geq 0$ for all $j\in \bigcup_{s\in y}s $, and we arrive at \eqref{infity_wasser_re_gamma}.
		
		%
		%
		%

		\item According to \eqref{infity_wasser_re_gamma}, we note that
		\begin{align*}
		v_{U\Gamma}&\leq \frac{1}{N}\sum_{k\in [N]}\max_{y\in F_\Gamma}\min_{s\in y}\sum_{j\in s}\bar{c}_j^k+\frac{1}{N}\sum_{k\in [N]}\max_{y\in F_\Gamma}\max_{\sum_{j\in \bigcup_{s\in y}s }(\beta_j^k)^r\leq \theta^r, \beta_j^k\geq 0,\forall j\in \bigcup_{s\in y}s }\max_{s\in y}\sum_{j\in s}\beta_j^k\\
		&=v_{U\Gamma}^{SAA}+\frac{1}{N}\sum_{k\in [N]}\max_{y\in F_\Gamma}\max_{s\in y}\max_{\sum_{j\in \bigcup_{s\in y}s }(\beta_j^k)^r\leq \theta^r, \beta_j^k\geq 0,\forall j\in \bigcup_{s\in y}s }\sum_{j\in s}\beta_j^k\\
		&=v_{U\Gamma}^{SAA}+\Gamma^{\frac{r-1}{r}}\theta
		\end{align*}
		where the first inequality is due to 
		\[\min_{s\in y}\sum_{j\in s}(\bar{c}_j^k+\beta_j^k)\leq \min_{s\in y}\sum_{j\in s}\bar{c}_j^k+\max_{s\in y}\sum_{j\in s}\beta_j^k,\]
		the first equality is due to the definition of $v_{U\Gamma}^{SAA}$ and the interchange of maximization operators, and the second one is due to the definition of the dual norm.
		
		On the other hand, replacing the inner minimization in \eqref{infity_wasser_re_gamma} by its lower bound, we have
			\begin{align*}
			v_{U\Gamma}&\geq \frac{1}{N}\sum_{k\in [N]}\max_{y\in F_\Gamma}\min_{s\in y}\sum_{j\in s}\bar{c}_j^k+\frac{1}{N}\sum_{k\in [N]} \min_{y\in F_\Gamma} \max_{\sum_{j\in \bigcup_{s\in y}s }(\beta_j^k)^r\leq \theta^r, \beta_j^k\geq 0,\forall j\in \bigcup_{s\in y}s }\min_{s\in y}\sum_{j\in s}\beta_j^k\\
			&\geq v_{U\Gamma}^{SAA}+\frac{1}{N}\sum_{k\in [N]}\min_{y\in F_\Gamma} \max_{\sum_{j\in \bigcup_{s\in y}s }(\beta_j^k)^r\leq \theta^r, \beta_j^k\geq 0,\forall j\in \bigcup_{s\in y}s }H(\bm{\beta})
			\end{align*}
			where $$H(\bm{\beta}):=\min_{\bm{\lambda}^k}\left\{\sum_{j\in \bigcup_{s\in y}s }\lambda_j^k\beta_j^k:\sum_{j\in \bigcup_{s\in y}s }\lambda_j^k=\Gamma,\lambda_j^k\in [0,1], \forall j\in \bigcup_{s\in y}s \right\}.$$
			Note that $H(\bm{\beta})$ is a symmetric concave function and the constraint system $\{\bm{\beta}:\sum_{j\in \bigcup_{s\in y}s }(\beta_j^k)^r\leq \theta^r, \beta_j^k\geq 0,\forall j\in \bigcup_{s\in y}s \}$ is also symmetric. Thus, the optimal value of
			\[\max_{\bm\beta }\left\{H(\bm{\beta}): \sum_{j\in \bigcup_{s\in y}s }(\beta_j^k)^r\leq \theta^r, \beta_j^k\geq 0,\forall j\in \bigcup_{s\in y}s\right\}\]
			is achieved by $\beta_j^k=\frac{\theta}{|\bigcup_{s\in y}s |^{\frac{1}{r}}}$ for each $j\in \bigcup_{s\in y}s $, i.e.,
			\begin{align*}
			v_{U\Gamma}&\geq v_{U\Gamma}^{SAA}+\frac{\Gamma\theta}{\max_{y\in F_{\Gamma}}|\bigcup_{s\in y}s |^{\frac{1}{r}}}.
			\end{align*}
		
		\QEDA
		
	\end{enumerate} 
	
	We remark that: (a) If there are a finite number of elements in $F_T$, then for each $y\in F_T$, we can convert the inner maximin problem in formulation \eqref{infity_wasser_re_gamma} as a convex program; and (b) In general, formulation \eqref{infity_wasser_re_gamma} can be difficult to solve. Thus, we recommend obtaining $v_{U\Gamma}^{SAA}$ instead and using inequalities \eqref{infity_wasser_re_gamma_bound} to estimate $v_{U\Gamma}$.
	
	The following results on confidence bounds follow the similar proof techniques in Theorem~\ref{thm_bound_u_gamma}, and thus their proofs are omitted for brevity.
	\begin{theorem}\label{thm_bound_u_gamma}Let $v_{U\Gamma}^{T}=\E_{\Pr^T}[\min_{x\in X}\max_{s\subseteq x, |s|=\Gamma}\sum_{j\in s}\tilde{c}_{j}]$ denote stochastic $\Gamma$BCP under the true distribution $\Pr^T$. Suppose that there exists a positive $\sigma$ such that $\E_{\Pr^T}[\exp((Z_{\Gamma}(\tilde{\bm c})-v_{U\Gamma}^{T})^2/\sigma^2)]\leq e$.
		Given $\epsilon\in (0,1)$, then
		\begin{enumerate}[(i)]
			\item let $\theta=N^{-\frac{1}{2}}\sigma \sqrt{-3\log (\epsilon)}\Gamma^{-1}{\max}_{y\in F_{\Gamma}}|\bigcup_{s\in y}s |^{\frac{1}{r}}=O(N^{-\frac{1}{2}})$. Then we have 
			$$\Pr^T\left\{v_{U\Gamma}\geq v_{U\Gamma}^T\right\}\geq 1-\epsilon;$$
			\item let $\theta=N^{-\frac{1}{2}}\sigma \sqrt{-3\log (\epsilon)}\Gamma^{-\frac{r-1}{r}}=O(N^{-\frac{1}{2}})$. Then we have 
			$$\Pr^T\left\{v_{U\Gamma}\leq v_{U\Gamma}^T+2\theta\right\}\geq 1-\epsilon.$$
		\end{enumerate}

	\end{theorem}
	
	\subsection{DR$\Gamma$BCP-D \eqref{general_eq_dis_robust_gamma_d}: Equivalent Reformulations, Confidence Bounds, and Alternative Decision-Robust Model}
	
	In this subsection, we will provide equivalent reformulations, and confidence bounds of DR$\Gamma$BCP-D \eqref{general_eq_dis_robust_gamma_d}. Similarly, DR$\Gamma$BCP-D \eqref{general_eq_dis_robust_gamma_d} may not yield robust decisions when the sample size is small, and therefore, we propose an alternative decision-robust model. The proofs in this subsection are almost identical to those in Section~\ref{sec_DRBCP-D} and are thus omitted.
	
	We first summarize the reformulation results.
	\begin{subequations}
		\begin{theorem}[Reformulations]\label{thm_drbcp_d_ref_gamma_d} Let $v_{D\Gamma}^{SAA}=\min_{x\in X}\frac{1}{N}\sum_{k\in[N]}\max_{s\subseteq x, |s|=\Gamma}\sum_{j\in s}\bar{c}_{j}^k$ denote sampling average approximation of stochastic decision-making $\Gamma$BCP with empirical distribution $\Pr_{\tilde{\bzeta}}$.
			Suppose $\|\cdot\|=\|\cdot\|_r$ with $r\geq 1$, then
			\begin{enumerate}[(i)]
				\item DR$\Gamma$BCP-D \eqref{general_eq_dis_robust_gamma_d} is equivalent to
				\begin{equation}\label{infity_wasser_re_gamma_d}
				v_{D\Gamma}=v_{D\Gamma}^{SAA}+\Gamma^{\frac{r-1}{r}}\theta
				\end{equation}
				\item Suppose that set $X$ admits a binary programming representation $\hat{X}\subseteq \{0,1\}^n$. Then $v_{D\Gamma}^{SAA}$ is equal to
				\begin{equation}\label{formulate_A_gamma}
				\begin{aligned}
				v_{D\Gamma}^{SAA}=\min_{\bm{z}^k\in \hat{X},\alpha^k,\bm{\pi}^k, \forall k\in [N]}&\frac{1}{N} \sum_{k\in [N]}\left(\alpha^k\Gamma+\sum_{j\in [n]}\pi_j^k\right)\\
				\text{s.t.}\quad&\alpha^k+\pi_j^k\geq (\bar{c}_j^k-L^k)z_j^k+L^k,\pi_j^k\geq 0, \forall j\in [n],k\in [N],
				\end{aligned}
				\end{equation}
				where $L^k:=\min_{\tau\in [n]}\bar{c}_\tau^k$.
				
			\end{enumerate} 
	\end{theorem}\end{subequations}
	
	Next, we summarize the confidence bounds.
	\begin{theorem}\label{thm_bound_d_gamma}Let $v_{D\Gamma}^{T}=\min_{x\in X}\E_{\Pr^T}[\sum_{k\in[N]}\max_{s\subseteq x, |s|=\Gamma}\sum_{j\in s}\tilde{c}_{j}]$ denote stochastic decision-making $\Gamma$BCP with true distribution $\Pr^T$. Suppose that there exists a positive $\sigma$ such that $$\E_{\Pr^T}\left[\exp\left(\sigma^{-2}\left(\max_{s\subseteq x, |s|=\Gamma}\sum_{j\in s}\tilde{c}_{j}-\E_{\Pr^T}\left[\sum_{k\in[N]}\max_{s\subseteq x, |s|=\Gamma}\sum_{j\in s}\tilde{c}_{j}\right]\right)^2\right)\right]\leq e$$ for each $x\in X$.
		Given $\epsilon\in (0,1)$ and $r\in [1,\infty]$, then
		\begin{enumerate}[(i)]
			\item let $\theta=\Gamma^{-\frac{r-1}{r}}N^{-\frac{1}{2}}\sigma \sqrt{-3\log (\epsilon)+3n\log(2)}=O(N^{-\frac{1}{2}})$. Then we have 
			$$\Pr^T\left\{v_{D\gamma}\geq v_{D\gamma}^T\right\}\geq 1-\epsilon;$$
			\item let $\theta=\Gamma^{-\frac{r-1}{r}}N^{-\frac{1}{2}}\sigma \sqrt{-3\log (\epsilon)+3n\log(2)}=O(N^{-\frac{1}{2}})$. Then we have 
			$$\Pr^T\left\{v_{D\gamma}\leq v_{D\gamma}^T+2\theta\right\}\geq 1-\epsilon.$$
		\end{enumerate}

	\end{theorem}
	
We propose the following decision-robust model as an alternative of DR$\Gamma$BCP-D \eqref{general_eq_dis_robust_gamma_d}.
	\begin{equation} \label{eq_var_for_0_2_gamma}
	\begin{aligned}
	v_{DR}^\gamma=\min_{x\in  \L_{\theta\gamma}(X)}\Var_{\Pr_{\tilde{\bzeta}}} \left[\max_{s\subseteq x, |s|=\Gamma}\sum_{j\in s}\tilde{c}_{j}\right],
	\end{aligned}
	\end{equation}
	where
	\begin{align}
	\L_{\theta\Gamma}(X)=\left\{x\in X:\frac{1}{N}\sum_{k\in [N]}\max_{s\subseteq x, |s|=\Gamma}\sum_{j\in s}\bar{c}_j^k\leq v_{D\Gamma}:=v_{D\Gamma}^{SAA}+\Gamma^{\frac{r-1}{r}}\theta \right\}.\label{eq_def_L_2}
	\end{align}
	
	Finally, we show that the decision-robust model \eqref{eq_var_for_0_2_gamma} can be formulated as an MISOCP.
	\begin{proposition}Suppose that set $X$ admits a binary programming representation $\hat{X}$. Then the decision-robust model \eqref{eq_var_for_0_2} is equivalent to the following MISOCP
		\begin{equation} \label{var_for_gamma}\small
		\begin{aligned}
		v_{DR}^\Gamma=\min_{\bm{z}\in \hat{X},\bm{v},\bm{\alpha},\bm{\pi},\bm{\lambda}}\ &\frac{1}{N}\sum_{k\in [N]}\left(v^k-\frac{1}{N}\sum_{\tau\in [N]}v^\tau\right)^2\\
		\text{s.t.}\quad& \frac{1}{N}\sum_{k\in [N]}v^k\leq v_{D\Gamma}^{SAA}+\Gamma^{\frac{r-1}{r}}\theta,\\
		& v^k=\alpha^k\Gamma+\sum_{j\in [n]}\pi_j^k \geq (\overline{c}_j^k-L^k)z_j+L^k,\forall \  j \in [n], k \ \in [N],\\
		& v^k \leq \sum_{j\in [n]} \lambda_j^k\overline{c}_j^k, \forall  k \in [N],\\
		& \sum_{j\in [n]} \lambda_j^k =\Gamma, \forall \  k \in [N],\\
		& \lambda_j^k \leq z_j, \forall \  k \in [N], j \ \in [n], \\
		& \lambda_j^k \in [0,1],  \forall \  j \in [n], k \ \in [N],
		\end{aligned}
		\end{equation}
		where $L^k:=\min_{\tau\in [n]}\bar{c}_\tau^k$.
	\end{proposition}
	%
	
	\section{Conclusion}\label{sec_conclusion}
	
	This paper studies distributionally robust bottleneck combinatorial problems (DRCBP) under the Wasserstein ambiguity set. Motivated by different applications, this paper focuses on two variants of DRCBP: DRCBP-U and DRCBP-D: (i) Motivated by the multi-hop wireless network, for DRCBP-U, the decision-makers would like to have a reliable estimation of the values of the stochastic bottleneck combinatorial problem; and (ii) Motivated by the ride-sharing transportation problem, for DRCBP-D, the decision-makers seek the best decision to minimize the worst-case bottleneck cost. In both models, the underlying true probability distributions are assumed to be not fully specified and are contained in an ambiguity set defined by a ball centered at empirical distribution within a certain Wasserstein distance. For both DRCBP-U and DRCBP-D, we derive their equivalent reformulations, demonstrate that they can be computed effectively as mixed-integer programs, and show that the Wasserstein radius can be chosen less conservatively by establishing their connections to the sampling average approximation counterparts. We show that the results can be extended to distributionally robust $\Gamma-$sum bottleneck combinatorial problems (DR$\Gamma$CBP), where the decision-makers would like to know the worst-case sum of $\Gamma$ largest elements instead of the bottleneck element.

{ 
	\section*{Acknowledgments}
The authors would like to thank Editor-in-Chief Dr. Sven Leyffer and two anonymous reviewers
for their valuable comments on improving this paper.
}
	
	\bibliography{bottleneck}

	\newpage
	\titleformat{\section}{\large\bfseries}{\appendixname~\thesection .}{0.5em}{}
	
	\begin{appendices}
		
		\section{DRBCP-U \eqref{general_eq_dis_robust} under $q-$Wasserstein Ambiguity Set: Equivalent Reformulation, MIP Representation, and Confidence Bounds}\label{1_norm_resul}
		For the completeness of this paper, in this appendix, we provide parallel results for DRBCP-U \eqref{general_eq_dis_robust} under $q-$Wasserstein ambiguity set \eqref{eq_general_das} with $q\in [1,\infty)$.

		\subsection{Equivalent Reformulation}
		We first derive a deterministic representation of DRBCP-U \eqref{general_eq_dis_robust} under $q-$Wasserstein ambiguity set.
		\begin{theorem}\label{thm_drbcp_d_ref_q} Suppose $\|\cdot\|=\|\cdot\|_r$ with $r\geq 1$, then DRBCP-U \eqref{general_eq_dis_robust} $q-$Wasserstein ambiguity set is equivalent to
			\begin{equation}\label{infity_wasser_re_q}
			v_U=\min_{\lambda\geq 0}\left\{\lambda\theta^q+\frac{1}{N}\sum_{k\in [N]}\max_{y\in F}\max_{t^k}\left\{t^k-\lambda \left(\sum_{j\in I^k(t^k,y)}(t^k-\bar c_j^k)^r\right)^{\frac{q}{r}}\right\}\right\},
			\end{equation}
			where $F$ is defined in Theorem \ref{cluter_blocker} and set $I^k(t,y)=\{j\in y: t-\overline{c}_j^k\geq 0\}$.
		\end{theorem}
		\proof
		According to Theorem~\ref{cluter_blocker} and theorem 1 in \cite{gao2016distributionally}, we have
		\begin{align*} 
		v_U&=\sup_{\Pr\in \P_{q}}\E_{\P}\left[Z(\tilde{\bm{c}})\right] 
		=\min_{\lambda\geq 0}\left\{\lambda\theta^q+\frac{1}{N}\sum_{k\in [N]}\max_{\bm{c}^k}\max_{y\in F}\left(\min_{j\in y}c_j^k-\lambda \|\bm{c}^k-\bar{\bm{c}}^k\|_r^q\right)\right\}.
		\end{align*}
		By swapping the two maximization operators and linearizing the inner minimization operator, we further have
		\begin{align*} 
		v_U&
		=\min_{\lambda\geq 0}\left\{\lambda\theta^q+\frac{1}{N}\sum_{k\in [N]}\max_{y\in F}\max_{t^k\leq c_j^k, \forall j\in y}\left(t^k-\lambda \|\bm{c}^k-\bar{\bm{c}}^k\|_r^q\right)\right\}.
		\end{align*}
		To simplify the inner maximization, we let $\beta_j^k=c_j^k-t^k$ for all $j\in y$ and define set $I^k(t,y)=\{j\in y: t-\overline{c}_j^k\geq 0\}$, then we have 
		\begin{align} 
		v_U&
		=\min_{\lambda\geq 0}\left\{\lambda\theta^q+\frac{1}{N}\sum_{k\in [N]}\max_{y\in F}\max_{\bm{\beta}^k,t^k}\left\{\left(t^k-\lambda \|t^k\e+\bm\beta^k-\bar{\bm{c}}^k\|_r^q:\beta_j^k \geq 0, \forall j\in y\right)\right\}\right\}\notag\\
		&
		=\min_{\lambda\geq 0}\left\{\lambda\theta^q+\frac{1}{N}\sum_{k\in [N]}\max_{y\in F}\max_{t^k}\left\{t^k-\lambda \left(\sum_{j\in I^k(t^k,y)}(t^k-\bar c_j^k)^r\right)^{\frac{q}{r}}\right\}\right\},\label{eq_lambda}
		\end{align}
		where the second equality is due to the optimality condition that if $j\in y\setminus I^k(t^k,y) $, we must have $\beta_j^k=-t^k+\bar c_j^k>0$; otherwise, $\beta_j^k=0$.

		\QEDA

		\subsection{Mixed Integer Programming Representation}
		Note that in Theorem \ref{thm_drbcp_d_ref_q}, the outer minimization problem is a univariate convex program and thus can be solved efficiently via the bisection method. We then observe that for any given $\lambda$, the inner maximizations in DRBCP-U \eqref{infity_wasser_re_q} can be decomposed into $N$ separate MIPs.
		\begin{proposition}\label{prop_mip_u_q}
			Suppose that the blocker $F$ admits a binary programming representation $\hat{F}\subseteq \{0,1\}^n$. Then
			$v_U=\min_{\lambda\geq0}\{\lambda\theta^q+1/N\sum_{k\in[N]}v_U^k(\lambda)\}$ and $v_U^k(\lambda)$ is equivalent to
			\begin{align}\label{mixed_integer_form_u_q}
			v_U^k(\lambda)=\max_{\bm{z}^k\in \hat{F}, t^k, \bm{\beta}^k,\theta^k}&\left\{t^k-\lambda (\theta^{k})^q: \|t^k\e-\bar {\bm c}_j^k+{\bm \beta}^k\|_r\leq \theta^{k},\beta_j^k\geq -M_j^k(1-z_j^k), \forall j\in [n]\right\},
			\end{align}
			where $M_{j}^k= \hat{M}^k-\bar{c}_{j}^k$ and
			$ \hat{M}^k=\max_{\tau\in [n]}\bar{c}_\tau^k+(\lambda q)^{-1/(q-1)}$.
		\end{proposition}
		\proof According to the proof of Theorem \ref{thm_drbcp_d_ref_q}, 	$v_U=\min_{\lambda\in [0,1]}\{\lambda\theta^q+1/N\sum_{k\in[N]}v_U^k(\lambda)\}$, where
		\[v_U^k(\lambda):=\max_{y\in F,t^k,\bm{\beta}^k,\theta^k}\left\{t^k-\lambda (\theta^{k})^q:\|t^k\e-\bar {\bm c}_j^k+{\bm \beta}^k\|_r\leq \theta^{k}, \beta_j^k\geq 0, \forall j\in y\right\}.\]
		Above, we can augment the vector $\bm{\beta}^k$ by defining free variable $\beta_j^k$ for each $j\in [n]\setminus y$ and rewrite $v_U^k(\lambda)$ as
		\[v_U^k(\lambda):=\max_{y\in F,t^k,\bm{\beta}^k,\theta^k}\left\{t^k-\lambda (\theta^{k})^q:\|t^k\e-\bar {\bm c}_j^k+{\bm \beta}^k\|_r\leq \theta^{k}, \beta_j^k\geq 0, \forall j\in y, \beta_j^k \in \Re, \forall  j\in [n]\setminus y\right\},\]
		since at optimality, we must have $\beta_{j}^k=-t^k+\bar c_j^k$ for each $j\in [n]\setminus y$.
		
		As the blocker $F$ admits a binary programming representation $\hat{F}$, the value $v_U^k(\lambda)$ is further equal to
		\begin{align}
		v_U^k(\lambda):=\max_{\bm{z}^k\in \hat{F},t^k,\bm{\beta}^k,\theta^k}\left\{t^k-\lambda (\theta^{k})^q:\|t^k\e-\bar {\bm c}_j^k+{\bm \beta}^k\|_r\leq \theta^{k}, z_j^k\beta_j^k\geq 0, \forall j\in [n]\right\}.\label{mixed_integer_form_u1_q}
		\end{align}
		To prove that \eqref{mixed_integer_form_u1_q} and \eqref{mixed_integer_form_u_q} are equivalent, we only need to show that $\beta_{j}^k\geq -M_j^k$ for each $j\in[n]$. Since at optimality, $\beta_{j}^k$ is equal to $0$ or $-t^k+\bar c_j^k$. Thus, it is sufficient to find an upper bound of $t^k$. In \eqref{infity_wasser_re_q}, suppose the optimal $t^k\geq\max_{\tau\in [n]}\bar{c}_\tau^k$, and then $I^k(t,y)=y$. The first order optimality condition of inner maximization of \eqref{infity_wasser_re_q} yields
		\begin{align*}
		1&= \lambda q\left(\sum_{j\in y}(t^k-\bar c_j^k)^r\right)^{\frac{q-r}{r}}\left(\sum_{j\in y}(t^k-\bar c_j^k)^{r-1}\right)\geq \lambda q\left(t^k-\max_{\tau\in [n]}\bar{c}_\tau^k\right)^{q-1},
		\end{align*}
		where the first inequality is due to $t^k-\bar c_j^k\geq t^k-\max_{\tau\in [n]}\bar{c}_\tau^k$ and $|y|\geq 1$. Thus, we can define $\hat{M}^k$ as
		\[ \hat{M}^k=\max_{\tau\in [n]}\bar{c}_\tau^k+(\lambda q)^{-1/(q-1)}.\]
		\QEDA
		
		We remark that: (i) to obtain $v_U$, one needs to first fix $\lambda$ and then solve $N$ separate MIPs \eqref{mixed_integer_form_u_q}, and then optimize $\lambda$, which requires more computational time compared to that of DRBCP-U \eqref{general_eq_dis_robust} under $\infty-$ Wasserstein ambiguity set; and (ii) as long as $r$ is a rational number, using the conic quadratic representation results in \cite{ben2001lectures}, then the MIP \eqref{mixed_integer_form_u_q} can be formulated as a mixed integer second order conic program (MISOCP).

		%
		%

		\subsection{Confidence Bounds}
		In this subsection, we plan to compare $v_U$ with its sampling average approximation counterpart, which further motivates us a choice of Wasserstein radius $\theta$. Recall that $v_U^{SAA}$ and $v_U^{T}$ are defined in \eqref{eq_saa} and \eqref{eq_saa_t}, respectively. Our goal is to determine a proper $\theta$ such that $v_U\geq v_U^{T}$ with a high probability. To begin with, we would like to show the relationship between $v_U$ and $v_U^{SAA}$.
		
		\begin{proposition} \label{infity_norm_saa_q} Let $v_U^{SAA}$ be defined in \eqref{eq_saa}. Then
			$$\frac{\theta}{\max_{y\in F}|y|^{\frac{1}{r}}}\leq v_U-v_U^{SAA}\leq \theta.$$
		\end{proposition}
		\proof 
		\begin{enumerate} [(i)]
			\item To prove $v^U-v^{SAA}\leq \theta$, 
			using the fact that $t^k \geq \min_{j\in y}\bar{c}_j^k$, then we can upper bound $v^U$ as
			\begin{align}
			v_U&\leq\min_{\lambda\geq 0}\left\{\lambda\theta^q+\frac{1}{N}\sum_{k\in [N]}\max_{y\in F}\max_{t^k\geq \min_{j\in y}\bar{c}_j^k}\left\{t^k-\lambda\left(t^k-\min_{j\in y}\bar{c}_j^k\right)^{q}\right\}\right\}.\label{eq_ub1}
			\end{align}
			To optimize the inner maximization, there are two cases:
			\begin{enumerate}[{Case} 1.]
				\item If $q=1$, then the inequality \eqref{eq_ub1} is equivalent to
				\begin{align*}
				v_U&\leq\min_{\lambda\geq 0}\left\{\lambda\theta+\frac{1}{N}\sum_{k\in [N]}\max_{y\in F} \min_{j\in y}\bar{c}_j^k: \lambda\geq 1\right\}=v^{SAA}+\theta.
				\end{align*}
				
				\item If $q>1$, then the inequality \eqref{eq_ub1} is equivalent to
				\begin{align*}
				v_U&\leq\min_{\lambda\geq 0}\left\{\lambda\theta^q+\frac{q-1}{q^{\frac{q}{q-1}}}\lambda^{-\frac{1}{q-1}}+\frac{1}{N}\sum_{k\in [N]}\max_{y\in F} \min_{j\in y}\bar{c}_j^k\right\}\\
				&=v^{SAA}+\theta.
				\end{align*}
				
			\end{enumerate}

			\item To prove $v^U-v^{SAA}\geq \frac{\theta}{{\max}_{y\in F}|y|^{\frac{1}{r}}}$, using the fact that $t^k \geq \min_{j\in y}\bar{c}_j^k$ again, then we can bound $v^U$ from the below as
			\begin{align}
			v_U&\geq\min_{\lambda\geq 0}\left\{\lambda\theta^q+\frac{1}{N}\sum_{k\in [N]}\max_{y\in F}\max_{t^k\geq \min_{j\in y}\bar{c}_j^k}\left\{t^k-\lambda|y|^{\frac{q}{r}}\left(t^k-\min_{j\in y}\bar{c}_j^k\right)^{q}\right\}\right\}.\label{eq_ub}
			\end{align}
			To optimize the inner maximization, there are two cases:
			\begin{enumerate}[{Case} 1.]
				\item If $q=1$, then the inequality \eqref{eq_ub} is equivalent to
				\begin{align*}
				v_U&\geq\min_{\lambda\geq 0}\left\{\lambda\theta+\frac{1}{N}\sum_{k\in [N]}\max_{y\in F} \min_{j\in y}\bar{c}_j^k: \lambda\geq \frac{1}{\min_{h\in F}|h|^{\frac{q}{r}}}\right\}=v^{SAA}+\frac{\theta}{\min_{y\in F}|y|^{\frac{q}{r}}}\\
				&\geq v^{SAA}+\frac{\theta}{\max_{y\in F}|y|^{\frac{q}{r}}}.
				\end{align*}
				
				\item If $q>1$, using the fact that $|y|^{\frac{q}{r}}\leq \max_{h\in F}|h|^{\frac{q}{r}}$ for all $y\in F$, then the right-hand side of inequality \eqref{eq_ub} is further lower bounded by
				\begin{align*}
				v_U&\geq\min_{\lambda\geq 0}\left\{\lambda\theta^q+\frac{q-1}{q^{\frac{q}{q-1}}}\left(\frac{1}{\lambda\max_{h\in F}|h|^{\frac{q}{r}}}\right)^{\frac{1}{q-1}}+\frac{1}{N}\sum_{k\in [N]}\max_{y\in F} \min_{j\in y}\bar{c}_j^k\right\}\\
				&=v^{SAA}+\frac{\theta}{\max_{y\in F}|y|^{\frac{1}{r}}}.
				\end{align*}
				\QEDA
			\end{enumerate}		
		\end{enumerate}
		
		%
		%
		%
		%
		%

		Now we are ready to show the following main result on the relationship between the value of DRBCP-U $v_U$ and true value $v_U^T$.
		\begin{theorem}\label{thm_bound_u_q}Suppose that there exists a positive $\sigma$ such that $\E_{\Pr^T}[\exp((Z(\tilde{\bm c})-v_U^T)^2/\sigma^2)]\leq e$.
			Given $\epsilon\in (0,1)$, then
			\begin{enumerate}[(i)]
				\item let $\theta=N^{-\frac{1}{2}}\sigma q^{-1/q}\sqrt{-3\log (\epsilon)}{\max}_{y\in F}|y|^{\frac{1}{r}}=O(N^{-\frac{1}{2}})$. Then we have 
				$$\Pr^T\left\{v_U\geq v_U^T\right\}\geq 1-\epsilon;$$
				\item let $\theta=N^{-\frac{1}{2}}\sigma q^{-1/q}\sqrt{-3\log (\epsilon)}=O(N^{-\frac{1}{2}})$. Then we have 
				$$\Pr^T\left\{v_U\leq v_U^T+2\theta\right\}\geq 1-\epsilon.$$
			\end{enumerate}

		\end{theorem}
		\proof The proof is almost identical to Theorem~\ref{thm_bound_u} and is thus omitted.
		\QEDA 
		
		
		
		
		\subsection{Numerical Illustration  of  DRBCP-U \eqref{general_eq_dis_robust} under $2-$Wasserstein Ambiguity Set}\label{appendix_num_2w}
		
		\begin{table}[htbp]
			\begin{center}
				\caption{Numerical Results of  DRBCP-U \eqref{general_eq_dis_robust} under $2-$Wasserstein Ambiguity Set with Application to the Multi-hop Network.} 
				\label{dro_value_w1}
				\small\setlength{\tabcolsep}{3.0pt}
				\begin{tabular}{c|c|cccc|c|c}
					\hline
					\multirow{2}{*}{$N$}  & \multirow{2}{*}{$\theta$} & \multicolumn{4}{c|}{$r=1$}                            & \multirow{2}{*}{$95\%$ CI of $v_U^{SAA}$}& \multirow{2}{*}{{Theoretical $95\%$ CI}} \\ \cline{3-6}
					&                           & Time     & $v_U$ & $\lambda^*$& $\theta^*$             &                                           \\ \hline
					\multirow{11}{*}{100} & 0.00             & 565.993  & 5.411 & 9761.176  & \multirow{11}{*}{0.14} & \multirow{11}{*}{{[}5.279, 5.543{]}}  & \multirow{11}{*}{{{[}5.375, 5.565{]}}}      \\
					& 0.02                      & 477.811  & 5.391 & 24.607    &                        &                                           \\
					& 0.04                      & 610.274  & 5.371 & 12.077    &                        &                                           \\
					& 0.06                      & 506.139  & 5.352 & 7.993     &                        &                                           \\
					& 0.08                      & 607.456  & 5.333 & 5.892     &                        &                                           \\
					& 0.10                      & 478.448  & 5.314 & 4.660     &                        &                                           \\
					& 0.12                      & 446.110  & 5.296 & 3.843     &                        &                                           \\
					& 0.14                      & 430.717  & 5.278 & 3.262     &                        &                                           \\
					& 0.16                      & 431.357  & 5.259 & 2.823     &                        &                                           \\
					& 0.18                      & 615.299  & 5.241 & 2.475     &                        &                                           \\
					& 0.20                      & 450.941  & 5.224 & 2.196     &                        &                                           \\ \hhline{========}
					\multirow{11}{*}{250} & 0.00                      & 1482.088 & 5.417 & 6180.340  & \multirow{11}{*}{0.10} & \multirow{11}{*}{{[}5.335, 5.499{]}}  & \multirow{11}{*}{{{[}5.394, 5.514{]}}}    \\
					& 0.02                      & 1617.128 & 5.398 & 24.728    &                        &                                           \\
					& 0.04                      & 1549.721 & 5.378 &12.235    &                        &                                           \\
					& 0.06                      & 1612.423 & 5.359 & 8.012     &                        &                                           \\
					& 0.08                      & 1725.868 & 5.339 & 6.002     &                        &                                           \\
					& 0.10                      & 1700.600 & 5.320 & 4.668     &                        &                                           \\
					& 0.12                      & 1715.304 & 5.302 & 3.857     &                        &                                           \\
					& 0.14                      & 1571.976 & 5.283 & 3.272     &                        &                                           \\
					& 0.16                      & 1571.375 & 5.264 & 2.823    &                        &                                           \\
					& 0.18                      & 1778.676 & 5.246 & 2.491     &                        &                                           \\
					& 0.20                      & 1734.860 & 5.228 & 2.233     &                        &                                           \\ \hhline{========}
					\multirow{11}{*}{500}    & 0.00                      & 2924.362 & 5.452 & 6180.340  & \multirow{11}{*}{0.06} & \multirow{11}{*}{{[}5.397, 5.508{]}} & \multirow{11}{*}{{{[}5.401, 5.486{]}}}      \\
					& 0.02                      & 3012.862 & 5.432 & 24.879    &                        &                                           \\
					& 0.04                      & 3492.183 & 5.413 & 12.407    &                        &                                           \\
					& 0.06                      & 3170.741 & 5.393 & 8.026     &                        &                                           \\
					& 0.08                      & 3496.983 & 5.374 & 6.002     &                        &                                           \\
					& 0.10                      & 3423.280 & 5.355 & 4.677     &                        &                                           \\
					& 0.12                      & 3416.526 & 5.336 & 3.897     &                        &                                           \\
					& 0.14                      & 3187.321 & 5.318 & 3.287     &                        &                                           \\
					& 0.16                      & 3284.556 & 5.300 & 2.823     &                        &                                           \\
					& 0.18                      & 3538.916 & 5.281 & 2.509     &                        &                                           \\
					& 0.20                      & 2845.976 & 5.263 & 2.255     &                        &                                           \\ \hline
				\end{tabular}
			\end{center}
		\end{table}
		
		\subsection{Numerical Illustration  of DRBCP-D \eqref{general_eq_dis_robust_d} under $\phi$-divergence}\label{appendix_num_phi}
		
		Alternatively, we can consider DRBCP-D \eqref{general_eq_dis_robust_d} under a $\phi$-divergence ambiguity set $\P$. In particular, we construct the ambiguity set based on the total-variational distance with radius $d\in [0,2]$. According to \cite{jiang2018risk}, given $d\in [0,2]$, the distributionally robust fair matching problem in Example \ref{example6} under total variation based ambiguity set admits the following form
			\begin{equation}\label{bap_div}
			\begin{aligned}
			v_{D\phi}=\min_{\bm{v,\bar{v}, z},\beta,r}\ &(1-d/2)\beta+\frac{1}{N}\sum_{k\in [N]}\bar v^k+dr/2\\
			\text{s.t.}\ & \bar v^k\geq v^k-\beta, \bar v^k\geq 0, \forall k \in [N]\\
			&r\geq v^k, \forall k \in [N], \\
			& v^k\geq \overline{c}_{ij}^kz_{ij}, \sum_{i\in [m]}z_{ij}=1, \sum_{j\in [m]}z_{ij}=1,  z_{ij}\in\{0,1\}, \  \forall i, j\in[m],\forall k\in [N], \\
			&z_{ij}\in\{0,1\} \  \forall i,j\in[m], \forall k \in [N].
			\end{aligned}
			\end{equation}
			
\begin{figure}[htbp]
				\centering
				\subfloat[0.48\textwidth][Mean v.s. $d$]{
					\includegraphics[width=0.48\textwidth]{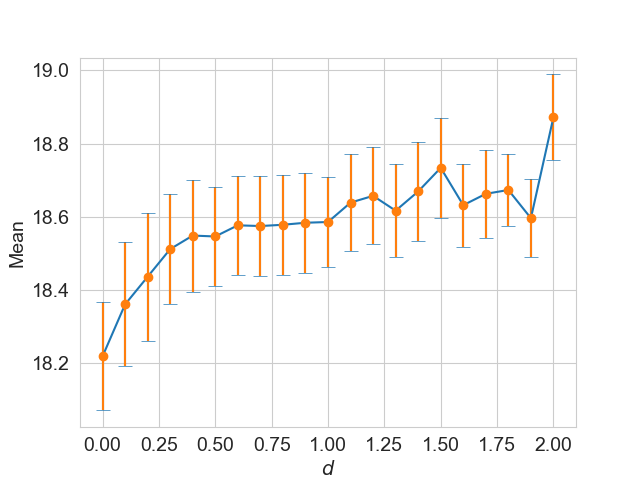} 
					\label{mean_50_div}
				}\hfill
				\subfloat[0.48\textwidth][Variance v.s. $d$]{
					\centering
					\includegraphics[width=0.48\textwidth]{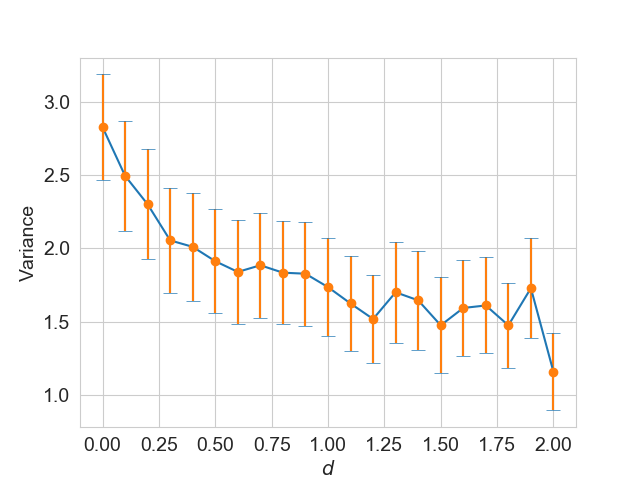} 
					\label{var_50_div}}
				\caption{Numerical Illustration of Fair Matching Problem under Total Variation-based Ambiguity Set \eqref{bap_div}: Using Real-world Dataset in Section \ref{sec_num_decision}}\label{fig_num_2_div}
			\end{figure}

				The numerical results using the same datasets as those in Section~\ref{sec_num_decision} are displayed in Figure \ref{fig_num_2_div} and Table~\ref{uber_match_div}.

	\begin{table}[htbp]
	\begin{center}
		\caption{Numerical Results of Fair Matching Problem under Total Variation-based Ambiguity Set \eqref{bap_div}: Using Hypothetical Dataset in Section \ref{sec_num_decision}}
		\label{uber_match_div}
		\small\setlength{\tabcolsep}{3.0pt}
		\small\setlength{\tabcolsep}{3.0pt}
		{ 
			\begin{tabular}{c|c|c|c|c|c|c|c|c|c|c}
				\hline
				\multirow{2}{*}{$N$}    & \multicolumn{5}{c|}{Case 1: $\N\left(\bm\mu_{NY},  10^2\bm\sigma_{NY}\right)$}                                                   & \multicolumn{5}{c}{Case 2: $\N\left(\bm\mu_{NY},  10^4\bm\sigma_{NY}\right)$}                                                       \\ \cline{2-11} 
				& $d$ & Time     & $v_{D\phi}$                    & $95\%$ CI  of  Solutions   &Theoretical $95\%$ CI     & $d$ & Time    & $v_{D\phi}$                      & $95\%$ CI  of  Solutions  &    Theoretical $95\%$ CI         \\ \hline
				\multirow{11}{*}{10} & 0.0   & 0.678 & 40.800      & $66.646\pm 0.130$    &   
				$66.646\pm 0.220$            & 0.0   & 0.571 & 242.350   & $497.251\pm 1.382$ &  $497.251\pm 2.345$ \\ \cline{2-6} \cline{7-11}
				& 0.2 & 0.579  & 45.790     & \multirow{4}{*}{$68.351\pm 0.127$}  & \multirow{4}{*}{$68.351 \pm 0.216$} 
				& 0.2 & 0.363 & 297.520   & $478.010\pm 1.253$     &  $478.010\pm 2.127$           \\ \cline{2-4} \cline{7-11} 
				& 0.4 & 0.394 & 49.070     &                                    
				& & 0.4 & 0.375 & 337.720    & \multirow{4}{*}{$468.365\pm 1.240$}&\multirow{4}{*}{$468.365\pm 2.104$}  \\ \cline{2-4} \cline{7-9}
				& 0.6 & 0.501 & 52.330     &                                     
				& & 0.6 & 0.551 & 369.560    &                                     \\ \cline{2-4} \cline{7-9}
				& 0.8 & 0.735 & 55.180     &                                    
				& & 0.8 & 0.308 & 397.080    &                                     \\ \cline{2-6} \cline{7-9}
				& 1.0   & 0.350 & 57.340     & $63.478\pm 0.130$                    
				&$63.478\pm 0.220$& 1.0   & 0.310 & 417.970    &                                     \\ \cline{2-6} \cline{7-11}
				& 1.2 & 0.327 & 59.350     & \bm{$63.362\pm 0.131$}                 
				&$63.362\pm  0.222$& 1.2 & 0.332 & 434.730    & \bm{$467.371\pm 1.276$}  & $467.371\pm 2.166$                \\ \cline{2-6} \cline{7-11}
				& 1.4 & 0.339   & 61.200      & \multirow{2}{*}{$63.478\pm 0.130$}   
				&\multirow{2}{*}{$63.478\pm 0.220$}& 1.4 & 0.334 & 447.180    & \multirow{2}{*}{$468.365\pm 1.240$}  & \multirow{2}{*}{$468.365\pm 2.104$} \\ \cline{2-4} \cline{7-9}
				& 1.6 & 0.281 & 62.590     &                                    
				&& 1.6 & 0.264  & 453.730    &                                     \\ \cline{2-6} \cline{7-11}
				& 1.8 & 0.050& 63.400     & $66.635\pm 0.127$                   
				&$66.635\pm 0.215$ & 1.8 & 0.079  & 458.400     & $505.127\pm 1.355$    & $505.127\pm 2.300$           \\ \cline{2-6} \cline{7-11}
				& 2.0   & 0.027 & 63.400      & {{$63.362\pm 0.131$} }              
				&$63.362\pm 0.222$  & 2.0   & 0.113     & 458.400     & {{$467.371\pm 1.276$}   }  &$467.371\pm 2.166$            \\ \hline\hline
				\multirow{11}{*}{20} & 0.0   & 1.606 & 45.930     & \multirow{7}{*}{$64.258\pm 0.165$}  &  \multirow{7}{*}{$64.258\pm 0.280$} & 0.0   & 1.685 & 304.910    & \multirow{5}{*}{\bm{$427.099\pm 1.149$}}& \multirow{5}{*}{$427.099\pm 1.950$} \\ \cline{2-4} \cline{7-9}
				& 0.2 & 1.1906 & 49.520     &                                    
				&& 0.2 & 2.357& 352.965   &                                     \\ \cline{2-4} \cline{7-9}
				& 0.4 & 0.920  & 52.950     &                                     
				&& 0.4 & 0.952 & 393.750    &                                     \\ \cline{2-4} \cline{7-9}
				& 0.6 & 1.082  & 55.600      &                                     
				&& 0.6 & 1.324 & 428.775   &                                     \\ \cline{2-4} \cline{7-9}
				& 0.8 & 0.800 & 57.875    &                                    
				& & 0.8 & 0.913 & 461.800     &                                     \\ \cline{2-4} \cline{7-11} 
				& 1.0   & 1.018  & 60.035    &                                   
				&  & 1.0   & 0.784 & 487.835   & \multirow{3}{*}{$501.812\pm 1.662$}& \multirow{3}{*}{$501.812\pm 2.821$} \\ \cline{2-4} \cline{7-9}
				& 1.2 & 0.750& 61.890     &                                    
				& & 1.2 & 0.548& 502.480    &                                     \\ \cline{2-9}
				& 1.4 & 0.865 & 63.420     & \multirow{2}{*}{$61.118\pm 0.159$}  
				& \multirow{2}{*}{$61.118\pm 0.271$} & 1.4 & 0.491 & 513.220    &                                     \\ \cline{2-4} \cline{7-11} 
				& 1.6 & 0.574 & 64.470     &                                    
				& & 1.6 & 0.521& 519.565   & \multirow{3}{*}{$431.159\pm 1.127$}& \multirow{3}{*}{$431.159\pm 1.912$} \\ \cline{2-9}
				& 1.8 & 0.319  & 64.825    & \bm{$59.228\pm 0.118$}                  
				& $59.228\pm 0.201$& 1.8 & 0.348 & 522.805   &                                     \\ \cline{2-9}
				& 2.0  & 0.034& 64.900      & $61.118\pm 0.159$                
				&   $61.118\pm0.271$& 2.0   & 0.036  & 523.000      &                                     \\ \hline\hline
				\multirow{11}{*}{40} & 0.0   & 4.829  & 49.4825   & \multirow{6}{*}{$61.118\pm 0.159$}  & \multirow{6}{*}{$61.118\pm 0.271$}  &
				0.0   & 4.217& 342.862   & $412.923\pm 1.086$& $412.923\pm 1.843$                   \\ \cline{2-4} \cline{7-11} 
				& 0.2 & 4.427  & 55.518   &                                   
				&  & 0.2 & 2.897  & 405.390    & $420.407\pm 1.133$   & $420.407\pm 1.922$                \\ \cline{2-4} \cline{7-11} 
				& 0.4 & 3.860 & 61.025    &                                     
				&& 0.4 & 1.636  & 455.753  & \bm{$397.311\pm 1.007$}  & {$397.311\pm 2.493$}              \\ \cline{2-4} \cline{7-11} 
				& 0.6 & 2.379  & 65.990   &                                    
				& & 0.6 & 1.661  & 493.1175  & \multirow{4}{*}{$415.991\pm 1.084$} & \multirow{4}{*}{$415.991\pm 1.840$}\\ \cline{2-4} \cline{7-9}
				& 0.8 & 2.576 & 70.675    &                                   
				&  & 0.8 & 1.193 & 526.960    &                                     \\ \cline{2-4} \cline{7-9}
				& 1.0   & 1.604  & 75.030     &                                     
				&& 1.0   & 1.239 & 556.585   &                                     \\ \cline{2-9}
				& 1.2 & 2.050& 78.975    & \multirow{3}{*}{\bm{$56.372\pm  0.103$}} 
				& \multirow{3}{*}{$56.372\pm 0.174$} & 1.2 & 1.256  & 582.03 0   &                                     \\ \cline{2-4} \cline{7-11} 
				& 1.4 & 3.237 & 81.685    &                                    
				& & 1.4 & 1.189   & 600.633  & \multirow{3}{*}{$466.052\pm 1.262$} & \multirow{3}{*}{$466.052\pm 2.143$}\\ \cline{2-4} \cline{7-9}
				& 1.6 & 1.553 & 83.970     &                                  
				&   & 1.6 & 0.954 & 612.990  &                                     \\ \cline{2-9}
				& 1.8 & 1.189& 85.838  & $63.167\pm 0.130$                  
				&  $63.167\pm0.221$& 1.8 & 0.622 & 621.578 &                                     \\ \cline{2-11} 
				& 2.0   & 0.062& 86.300      & $60.244\pm 0.118$                  
				&$60.244 \pm 0.201$ & 2.0  & 0.044 & 626.500     & $429.502\pm 1.092$ & $429.502\pm 1.854$                  \\ \hline
			\end{tabular}
		}
	\end{center}
\end{table}
	
	\end{appendices}

\end{document}